\documentclass{article}

\usepackage{microtype}
\usepackage{graphicx}
\usepackage{subfigure}
\usepackage{booktabs} 

\usepackage{hyperref}
    \hypersetup{colorlinks=true, linkcolor=black, citecolor=black, urlcolor=black}

\usepackage{microtype}

\usepackage[shortlabels, inline]{enumitem} 

\usepackage[retainorgcmds]{IEEEtrantools}
\usepackage{multirow}

\usepackage[noend,ruled,linesnumbered]{algorithm2e}
\DontPrintSemicolon 
\SetKwInput{KwInput}{Input}
\SetKwInput{KwOutput}{Output}

\usepackage{setspace}

\usepackage[normalem]{ulem}

\usepackage{amssymb}
\usepackage{amsmath}
\usepackage{amsthm}
\usepackage{bbm} 
\usepackage{bbold}
\newcommand{\1}{\mathbb{1}}
\providecommand{\1}{\mathds{1}}
\newcommand{\0}{\mathbb{0}}

\usepackage{pgfplots}
\usepackage{pgfplotstable}
\usepackage{tikz}
    \definecolor{myblue}{HTML}{3333CC}
    \definecolor{mylightblue}{HTML}{CCCCFF}
    \definecolor{myorange}{HTML}{FF6600}
    \definecolor{myred}{HTML}{CC3333}
    \definecolor{mygreen}{HTML}{11AA11}
    \definecolor{mylightgray}{HTML}{E0E0E0}
    \usetikzlibrary{decorations.pathmorphing}
    \tikzstyle{vertex}=[circle, draw, fill=white, inner sep=0pt, minimum width=1ex]
    \tikzset{every picture/.append style={baseline,scale=1.1}}
    \tikzstyle{cut-edge}=[dotted]
    \tikzstyle{base}=[line width=1.2pt]
    \tikzstyle{lifted}=[blue]

\usetikzlibrary{patterns}
\usetikzlibrary{external}
\usetikzlibrary{arrows,arrows.meta,bending}

\providecommand{\myparagraph}[1]{\textbf{#1}\hspace{0.5em}}

\providecommand{\R}{\mathbb{R}}

\providecommand{\N}{\mathbb{N}}

\providecommand{\C}{\mathcal{C}}

\providecommand{\F}{\mathcal{F}}

\providecommand{\X}{\mathcal{X}}
\providecommand{\I}{\mathcal{I}}
\providecommand{\J}{\mathcal{J}}

\providecommand{\bddtrue}{\top}
\providecommand{\bddfalse}{\bot}

\providecommand{\la}{\langle}
\providecommand{\ra}{\rangle}

\providecommand{\abs}[1]{\lvert #1 \rvert}

\theoremstyle{plain}

\newtheorem{prop}{Proposition}

\theoremstyle{definition}
\newtheorem{defn}{Definition}

\newtheorem{ex}{Example}

\theoremstyle{remark}
\newtheorem*{rem}{Remark}

\DeclareMathOperator{\sign}{sign}
\DeclareMathOperator{\idx}{idx}

\providecommand{\setsep}{\mid}

\usepackage{twoopt}
\usetikzlibrary{fit}
\usetikzlibrary{tikzmark}

\newcommand\drawCodeBox[2]{%
\begin{tikzpicture}[remember picture,overlay]
\coordinate (start) at ([yshift=1.4ex]pic cs:#1);
\coordinate (end) at ([yshift=-0.5ex,xshift=2.3ex]pic cs:#2);
\node[inner sep=2pt,draw=black,fit=(start) (end)] {};
\end{tikzpicture}%
}

\usepackage{hyperref}

\usepackage{etoolbox}

\usepackage[accepted]{icml2021}

\icmltitlerunning{Efficient Message Passing for 0--1 ILPs with Binary Decision Diagrams}

\begin{document}

\twocolumn[
\icmltitle{Efficient Message Passing for 0--1 ILPs with Binary Decision Diagrams}



\icmlsetsymbol{equal}{*}

\begin{icmlauthorlist}
\icmlauthor{Jan-Hendrik Lange}{unitueb}
\icmlauthor{Paul Swoboda}{mpi}
\end{icmlauthorlist}

\icmlaffiliation{mpi}{Max Planck Institute for Informatics, Saarbr\"ucken, Germany}
\icmlaffiliation{unitueb}{University of T\"ubingen, Germany}

\icmlcorrespondingauthor{Paul Swoboda}{pswoboda@mpi-inf.mpg.de}

\icmlkeywords{}

\vskip 0.3in
]



\printAffiliationsAndNotice{}  


\begin{abstract}
We present a message passing method for 0--1 integer linear programs.
Our algorithm is based on a decomposition of the original problem into subproblems that are represented as binary decision diagrams.
The resulting Lagrangean dual is solved iteratively by a series of efficient block coordinate ascent steps.
Our method has linear iteration complexity in the size of the decomposition and can be effectively parallelized.
The characteristics of our approach are desirable towards solving ever larger problems arising in structured prediction.
We present experimental results on combinatorial problems from MAP inference for Markov Random Fields, quadratic assignment, discrete tomography and cell tracking for developmental biology and show promising performance.
\end{abstract}


\section{Introduction}
\label{sec:introduction}

Structured prediction tasks in machine learning commonly require solving relaxations of NP-hard combinatorial optimization problems that are formulated as integer linear programs (ILPs).
Examples include discrete graphical models~\cite{werner2007linear}, graph partitioning~\cite{bansal2004correlation}, graph matchings~\cite{torresani2008feature} and tracking problems~\cite{luo2014multiple}.
Commercial ILP solvers like Gurobi~\cite{gurobi} or CPLEX~\cite{cplex} rely on standard linear programming algorithms, such as the simplex and barrier method.
These methods require matrix factorization and hence have super-linear complexity, which diminishes their competitiveness for very large problems that arise in structured prediction.
Therefore, considerable research effort has been invested into efficient dedicated solvers for specific problem classes.
Some of the most scalable methods exploit problem decompositions in order to solve Lagrangean relaxations. 
This family of algorithms includes subgradient, Frank-Wolfe and dual block coordinate ascent (DBCA) methods.
The DBCA approach, also called \emph{message passing} in the literature, exhibits very good performance for certain classes of combinatorial problems such as inference in Markov Random Fields, significantly outperforming general-purpose ILP solvers~\cite{Kappes2015}.
The drawback of DBCA solvers, however, is that they are only applicable to their dedicated problem class.
In order to solve other problem classes, new specialized DBCA algorithms have to be developed.
This hinders application of structured prediction in machine learning, since algorithm design is challenging and implementation time-consuming.

In this work we propose an efficient message passing method for solving relaxations of 0--1 ILPs.
Our method
(i)~has linear iteration complexity unlike LP solvers,
(ii)~is not restricted to a narrow subclass of problems unlike specialized solvers and
(iii) can be effectively parallelized and shows significant parallelization speedups.
We demonstrate the potential of our method on a wide variety of structured prediction tasks.

Our algorithm works by decomposing any given 0--1 ILP into smaller subproblems represented by binary decision diagrams (BDDs)~\cite{bryant1986graph}.
In the basic version we generate a BDD for each row of the constraint matrix.
    While general linear constraints may lead to BDDs of intractable sizes (e.g.\ for the NP-hard Knapsack problem), many inequalities commonly encountered in practice admit tractable BDD representations \cite{knuth2011art,wegener2000branching}.
Inequalities that do not admit a small BDD representation can be represented as multiple BDDs~\cite{abio2012new} efficiently.
Hence, our BDD-representation can always be chosen bounded by the size of the original problem description.
We combine the BDD subproblems via Lagrangean dual variables in order to obtain a convex relaxation for the ILP.
The algorithm then updates dual variables iteratively by an operation called min-marginal averaging, which maximally improves the objective w.r.t.\ the current set of dual variables.
Hence, our method belongs to the family of DBCA algorithms.
Based on the computed dual solution, we compute a primal one via depth-first exploration of the search space.
Our primal heuristic uses dual costs to guide the search towards high-quality solutions and individual subproblems inform the search so that feasible solutions are found fast and are of high quality.
We show how BDDs support an efficient implementation of our DBCA method.

We present all proofs in the appendix.
The code and datasets are available on \url{https://github.com/LPMP/BDD}.


\section{Related Work}
\label{sec:related-work}

\myparagraph{Dual Block-Coordinate Ascent} 
In machine learning DBCA algorithms for Lagrangean decompositions of combinatorial problems were successfully applied to a number of different tasks such as multiple target tracking~\cite{arora2013higher}, graph matching (quadratic assignment problem)~\cite{zhang2016pairwise,swoboda2017study}, multi-graph matching~\cite{swoboda2019convex}, the multicut problem~\cite{swoboda2017message}, cell tracking in biological image analysis~\cite{haller2020primal} and MAP inference in MRFs~\cite{kolmogorov2006convergent,kolmogorov2014new,globerson2008fixing,werner2007linear,savchynskyy12efficient,jancsary2011convergent,meltzer2012convergent,wang2013subproblem,johnson2007lagrangian,tourani2018mplp++,tourani2020taxonomy}.

The study \cite{swoboda2017dual} presents general algorithmic principles on how to design efficient DBCA algorithms for arbitrary combinatorial subproblems.
However, the decomposition into subproblems, the specific choice of updates and their efficient implementation are still left open for the algorithm designer to decide anew for each new problem class.
The work~\cite{werner2019relative} analyzes different update operations for DBCA problems and theoretically characterizes the resulting fixed points.

\myparagraph{Optimization with Binary Decision Diagrams} 
While binary decision diagrams (BDDs) have been used mostly to encode Boolean functions, finding optimal assignments w.r.t.\ a linear cost is a straightforward extension~\cite{knuth2011art}.
However, for NP-hard combinatorial problems the size of the BDD encoding  increases exponentially, which makes a straightforward application of BDDs computationally intractable.
In this context, two approaches have been proposed to limit BDD growth:
(i)~Solving a relaxation in which the set of feasible points is larger which leads to a lower bound~\cite{andersen2007constraint} or
(ii)~solving a restriction in which only a subset of feasible solutions is considered, leading to an upper bound~\cite{bergman2016decision,bergman2016discrete}.
In any case, the constructed BDD is significantly smaller than the BDD encoding of the original problem.

Related to our work is the Lagrangean relaxation method by \citet{bergman2015lagrangian}, who combine relaxed multi-valued decision diagrams (MDDs), an extension of BDDs, with additional constraints such as linear inequalities.
The resulting Lagrangean dual problem is solved with subgradient ascent.
In contrast, our method combines an arbitrary number of BDDs into a Lagrangean dual and applies the more efficient DBCA method.
\citet{improved_job_sequence_bounds_hooker_2019} extends the work of \citet{bergman2015lagrangian} to provide improved bounds on job sequencing problems.
Similarly, \citet{castro2020mdd} combine a relaxed decision diagram and linear constraints into a Lagrangean dual for a routing problem.
The integration of techniques based on decision diagrams and (mixed-)integer programming is considered by \cite{tjandraatmadja2020incorporating,gonzalez2020bdd,gonzalez2020integrated}. The latter apply their hybrid approach to the (quadratic) stable set problem.
Similar to our work, \citet{BergmanCire2016,bergman2018discrete,lozano2018consistent} consider a decomposition based on a collection of BDDs.
Their approach is to derive a lifted linear formulation from the intersection of network flow polytopes associated with the individual BDDs solved subsequently by an ILP solver.

In contrast to prior BDD-based work, which is either applicable beyond narrow problem classes or standalone, the method we propose is both at once.
It is (in principle) applicable to general 0--1 ILPs and does not rely on any external solvers.

\myparagraph{Integer Linear Programming} 
The integer linear programming (ILP) approach has been pioneered for the traveling salesman problem~\cite{dantzig1954solution} and is successful in solving many large scale instances~\cite{applegate2006traveling}.
The input to ILP solvers consists of a description of the feasible set in terms of linear inequalities.
In the first step, a linear programming (LP) relaxation obtained from relaxing the integrality constraints is solved.
If the relaxed solution has fractional entries, then cutting plane methods seek to tighten the relaxation.
Additionally, branch-and-bound steps are employed, which (in a nutshell) fix subsets of variables to integer values and resolve the modified LP.
In this way, the search space of integer feasible solutions is traversed recursively.
Branches of the search tree can be discarded if their associated lower bounds exceed the best found primal solution value.

Over the years, ILP technology has made spectacular progress.
In combination with hardware improvements, current state-of-the-art solvers outperform earlier ones by at least 6 orders of magnitude.
On a machine-independent basis, \citet{Bixby2012} estimates a speedup of a factor of $29000$ in the timespan 1991--2007 for CPLEX~\cite{cplex} and of a factor of $16$ in the timespan 2009--2012 for Gurobi~\cite{gurobi}.
With further advancements since then, solvers are routinely capable of solving previously inaccessible instances in moderate time.
Due to the considerable implementation effort, solvers written in academic projects seem not able to achieve state-of-the-art results~\cite{mittelmann2017latest,Mittelmann2020,MittelmannBenchmark}.


In the context of very large scale instances, a major challenge for ILP solvers is their limited capability to utilize modern parallel CPU architectures.
Solving LP relaxations, which consumes a considerable portion of the overall computation time, is particularly difficult to parallelize.
Although there have been efforts to exploit parallelization in the (dual) simplex method \cite{huangfu2018parallelizing}, the performance improvement over leading sequential implementations is relatively small~\cite{gurobi}.
The barrier method can be more effectively parallelized~\cite{gondzio2003parallel}, but often falls short of the dual simplex method on large sparse problems.
Moreover, the dual simplex method is more suitable for solving relaxations of ILPs due to its advantages in reoptimization.
We show in the experimental section that our method benefits significantly from parallelization when solving large scale problems and scales favourably w.r.t.\ problem size.


\section{Dual Decomposition of 0--1 Programs}
\label{sec:decomposition}
Below we introduce a Lagrangean decomposition for binary programs and recapitulate the min-marginal averaging algorithm, the most commonly used DBCA approach.

\begin{defn}[Binary program]
Consider $m$ index subsets~$\I_j \subset [n] = \{1, \dotsc, n\}$ and corresponding constraints $\X_j \subset \{0,1\}^{\I_j}$ for $j \in [m]$ together with a linear objective $c \in \R^n$.
The corresponding binary program is defined as
\begin{align}
\min \quad & c^\top x \quad \text{s.t.} \quad x_{\I_j} \in \X_j \quad \forall j \in [m]. \label{eq:binary-program} \tag{BP} 
\end{align}
\end{defn}
While we focus on 0--1 integer linear programs, binary programs can also encode, e.g., Max-SAT and weighted constraint satisfaction problems with finite domains.
\begin{ex}
	\label{ex:ILP}
Consider the 0--1 integer linear program
\begin{align}
\min \quad & c^\top x \quad
\text{s.t.} \quad Ax \leq b, \; x \in \{0,1\}^n.
\tag{ILP}
\label{eq:ILP}
\end{align}
The system of linear constraints $Ax \leq b$ may be split into $m$ blocks, each block representing a single (or multiple) rows of the system. For instance, let $a_j^\top x \leq b_j$ denote the $j$-th row of $Ax \leq b$, then the problem can be written in the form \eqref{eq:binary-program} by setting $\I_j = \{ i \in [n] \setsep a_{ji} \neq 0\}$ and $\X_j = \{ x \in \{0,1\}^{\I_j} \setsep \sum_{i \in \I_j} a_{ji} x_{i} \leq b_j \}$.
\end{ex}


\subsection{Lagrangean Dual}

The problem \eqref{eq:binary-program} is NP-hard and thus difficult to solve in general.
However, optimization over a single constraint alone, i.e.\ $\min_{x \in \X_j} c^\top x$, may be much easier depending on the constraint (although still NP-hard in general).
We exploit efficient methods for the individual subproblems in order to solve a Lagrangean dual problem of \eqref{eq:binary-program} by block coordinate ascent.

\begin{defn}[Lagrangean dual problem]
Define the set of subproblems that constrain variable $x_i$ as
$\J_i = \{j \in [m] \mid i \in \I_j\}$.
Let the energy for subproblem $j \in [m]$ w.r.t.\ Lagrangean dual variables $\lambda^j \in \R^{\I_j}$ be
$E^j(\lambda^j) = \min_{x \in \X_j} x^\top \lambda^j$.
Then the Lagrangean dual problem is defined as
\begin{align}
\max_\lambda \quad & \sum_j E^j(\lambda^j) \quad
\text{s.t.} \quad \sum_{j \in \J_i} \lambda^j_i = c_i \quad \forall i \in [n]. \label{eq:dual-problem} \tag{D}
\end{align}
\end{defn}

If the optima of the individual subproblems $E^j(\lambda^j$) agree with each other, then the consensus vector solves the original problem~\eqref{eq:binary-program} due to the constraints on $\lambda$.
In general, it is a lower bound on~\eqref{eq:binary-program}, since such a consensus need not hold.
We provide a formal derivation of problem \eqref{eq:dual-problem} in the appendix.

\subsection{Min-Marginal Averaging}
\label{sec:message-passing}

In this section we present the block coordinate ascent method to solve problem~\eqref{eq:dual-problem}.
The underlying algorithmic idea is to iterate over the variable indices $i \in [n]$ and update the associated dual variables such that in each subproblem the minima w.r.t.\ the current $x_i$ agree with each other.
This results in an algorithm which produces a monotonically non-decreasing sequence of lower bounds for~\eqref{eq:dual-problem}.
To this end, we define min-marginals and the min-marginal averaging update step. 

\begin{defn}[Min-marginal averaging]
For $i \in [n]$, $j \in \J_i$ and $\beta \in \{0,1\}$ let
\begin{align}
m^\beta_{ij} = \min_{x \in \X_j} x^\top \lambda^j \quad \text{s.t.} \quad x_i = \beta \label{eq:min-marginals}
\end{align}
denote the \emph{min-marginal} w.r.t.\ $i,j$ and $\beta$.
The \emph{min-marginal averaging update} w.r.t.\ $i$ is defined as
\begin{align}
\lambda^j_i \leftarrow \lambda^j_i - (m^1_{ij} - m^0_{ij}) + \frac{1}{\abs{\J_i}} \sum_{k \in \J_i} m^1_{ik} - m^0_{ik} \label{eq:mma-update}
\end{align}
for all $j \in \J_i$.
\end{defn}
The quantity $\abs{m^1_{ij} - m^0_{ij}}$ indicates by how much $\min_{x \in \X_j} x^\top \lambda^j$ increases if $x_i$ is fixed to $1$ (if $m^1_{ij} > m^0_{ij}$), respectively $0$ (if $m^1_{ij} < m^0_{ij}$).
The min-marginal averaging update results in an equal distribution of min-marginal differences across each involved subproblem.
Also, it results in a maximal improvement of the lower bound given all other dual variables are fixed.

\begin{prop}
\label{prop:mma-update}
The min-marginal averaging update w.r.t.\ $i \in [n]$ increases the dual bound by the non-negative value
\begin{align*}
 \min \Big\{ 0, \sum_{j \in \J_i} m^1_{ij} - m^0_{ij} \Big\} - \sum_{j \in \J_i} \min \left\{0,  m^1_{ij} - m^0_{ij} \right\}.
\end{align*}
\end{prop}

\paragraph{Min-Marginal Averaging Algorithm}
Algorithm~\ref{alg:min-marginal-averaging} computes a solution of the Lagrangean dual~\eqref{eq:dual-problem} by performing a series of min-marginal update steps~\eqref{eq:mma-update}.
The order in which the variables are processed, i.e.\ an appropriate permutation of $[n]$, is fixed at the start of the algorithm.
A suitable order can be obtained, for instance, from bandwidth minimization of the constraint matrix, cf.\ Section~\ref{sec:variable-order} in the appendix.
In the appendix we also present an alternative averaging strategy (Section~\ref{sec:averaging-strategy}) and a smoothed variant of our method (Section~\ref{sec:smoothing}) that is valuable for problems with bad fixed points of the non-smooth algorithm.

\begin{algorithm}
\caption{Min-marginal averaging}
\label{alg:min-marginal-averaging} 
\textbf{input} objective vector $c \in \R^n$, constraint sets $\X_j \subset \{0,1\}^{\I_j}$ for $j \in [m]$ \;
Find variable ordering $\{i_1,\ldots,i_n\} = [n]$.\;
Initialize dual variables $\lambda^j_i = c_i / \abs{\J_i}$ for all $i \in [n]$ and $j \in \J_i$.\;
    \While{(stopping criterion not met)}{
	Perform forward pass:\;
    \For {$i = i_1,\ldots,i_n$}{
        \For {$j \in \J_i$}{
            \tikzmark{mm-start}%
            Compute min-marginals for $\beta \in \{0,1\}$:\label{alg:min-marginal-averaging-mma-update-start} \;
$m^\beta_{ij} = \min_{x \in \X_j} x^\top \lambda^j$\ \  s.\ t.\ \ \ $x_i = \beta$
			\label{alg:min-marginal-averaging-mma-update-end} \tikzmark{mm-end}\; 
        }
\drawCodeBox{mm-start}{mm-end}%
        \For {$j \in \J_i$}{
            Update dual variable: $\lambda^j_i \leftarrow \lambda^j_i - (m^1_{ij} - m^0_{ij}) + \frac{1}{\abs{\J_i}} \sum_{k \in \J_i} m^1_{ik} - m^0_{ik}$ .\;
        }
    }
	Perform backward pass analogously (set variable order to $\{i_n,\ldots,i_1\}$)
    }
\end{algorithm}


\section{Primal Heuristic}
\label{sec:variable-fixing}

In this section we present our approach to determine good primal solutions for the problem \eqref{eq:binary-program} based on the dual solution that we obtained by block coordinate ascent.
The basic idea is to iteratively fix primal variables and backtrack until we arrive at a feasible solution.
The success of the search is determined by the order and the values of the variable fixations.
In order to define a search strategy, we compute variable scores $S_i \in \R$ and preferred variable values $\beta_i \in \{0,1\}$ for all $i \in [n]$.
Variables are fixed to preferred values $\beta_i$ in descending order of their scores $S_i$.

We describe a straightforward choice for $S_i$ and $\beta_i$ here and discuss alternatives in the appendix.
To this end, we compute for every index $i \in [n]$ the total min-marginal difference defined below.

\begin{defn}[Total min-marginal difference]
The \emph{total min-marginal difference} for $i \in [n]$ w.r.t.\ current dual variables~$\lambda$ is defined as
$M_i = \sum_{j \in \J_i} m^1_{ij} - m^0_{ij}$. 
\end{defn}
The quantity $\abs{M_i}$ indicates by how much the dual bound increases if $x_i$ is fixed to $1$ (if $M_i > 0$), respectively $0$ (if $M_i < 0$), in total across all individual subproblems, given that the dual variables remain unchanged.
Thus, we prefer to fix $x_i$ as indicated by the sign of $M_i$ and define
\begin{align}
\label{eq:preferred-variable-value}
\beta_i = 1 \quad \text{if} \quad M_i \leq 0 \quad \text{and} \quad \beta_i = 0 \quad \text{if} \quad M_i > 0.
\end{align}
The two variable fixation orders we employ are based on setting $S_i = \abs{M_i}$ or $S_i = -M_i$.
The second, less intuitive order introduces a bias to fix variables to $1$, which, in our experiments, allows to find feasible solutions faster.
This is due to feasible solutions typically having a small number of $1$-entries.

Our primal heuristic is detailed in Algorithms~\ref{alg:primal-heuristic} and~\ref{alg:restriction-propagation} in the appendix.
Algorithm~\ref{alg:primal-heuristic} traverses the space of feasible solutions by depth-first search.
We accelerate the search further by propagation of the restrictions to the feasible sets given the current partial assignment, see Algorithm~\ref{alg:restriction-propagation}.


\section{Implementation with BDDs}
\label{sec:bdd}

We have described above a generic DBCA procedure for optimizing a Lagrangean relaxation of 0--1 integer linear programs and a primal search heuristic for finding solutions given dual variables.
Below, we describe how reduced ordered binary decision diagrams~\cite{bryant1986graph}, a data-structure for representing Boolean functions, can provide efficient procedures for all operations that we require in the algorithms above.
First we define the type of binary decision diagrams we employ.

\subsection{Binary Decision Diagrams}

\begin{defn}[Binary decision diagram]
\label{def:bdd}
Given a set of ordered variable indices $\I = \{i_1,\ldots,i_k\} \subseteq [n]$, a \emph{binary decision diagram} (BDD) is a rooted, directed acyclic graph $(V, A)$ with $A = A^1 \cup A^0$ such that:
	\begin{enumerate}[(i),noitemsep,topsep=0pt,parsep=0pt,partopsep=0pt,leftmargin=0.7cm]
\item Every node $v \in V$ is labeled with an associated index $\idx(v) \in \I$ or is one of the two special nodes $\bddtrue$ or $\bddfalse$, representing the outcomes true and false, respectively.
\item The root node $r$ has index $\idx(r) = i_1$.
\item Each node $v \in V \backslash \{\bddtrue,\bddfalse\}$ has two outgoing arcs, a $1$-arc $(v,v^+_1) \in A^1$ and a $0$-arc $(v,v^+_0) \in A^0$.
\item Every path from $r$ to $\top$ or $\bot$ traverses a sequence of nodes $(v_1,\ldots,v_\ell,\bddfalse)$ for $\ell \leq k$ or $(v_1,\ldots,v_k,\bddtrue)$ with consecutive indices $\idx(v_j) = i_j$ for all $j$.
\item From every node $v \in V \backslash \{\bddtrue,\bddfalse\}$ there exists a path from $r$ to $v$ and from $v$ to $\bddtrue$.
\item The BDD is \emph{reduced}, i.e.\ there are no isomorphic subgraphs.
\end{enumerate}
\end{defn}

\begin{figure}[!ht]
\center
\begin{tikzpicture}[xscale=1.7,yscale=1.5,xslant=0.0]
	\tikzstyle{vertex}=[xslant=0.0,draw,circle,inner sep=1pt, minimum width=0pt]         
        
	\node[vertex] (1) at (0,0.5) {$1$};
	\node[vertex] (3a) at (1,0) {$3$};
	\node[vertex] (3b) at (1,1) {$3$};
	\node[vertex] (7a) at (2,0) {$7$};
	\node[vertex] (7b) at (2,1) {$7$};
	
	\node[vertex] (false) at (3,0) {$\bot$};
	\node[vertex] (true) at (3,1) {$\top$};

	\draw (1) edge[very thick,mygreen,->,>=stealth,dotted] node[black,below left] {\small $x_1=0$} (3a);
	\draw (1) edge[->,>=stealth] node[above left] {\small $x_1=1$} (3b);
	\draw (3a) edge[->,>=stealth,dotted] (7a);
	\draw (3a) edge[very thick,mygreen,->,>=stealth] (7b);
	\draw (3b) edge[->,>=stealth,dotted] (7b);
	\draw (3b) edge[->,>=stealth,bend left=10] (false);
	\draw (7a) edge[->,>=stealth,dotted] (false);
	\draw (7a) edge[->,>=stealth] (true);
	\draw (7b) edge[very thick,mygreen,->,>=stealth,dotted] (true);
	\draw (7b) edge[,->,>=stealth,bend left=15] (false);
		
\end{tikzpicture}
\caption{
    Binary decision diagram representing the indicator function of all binary triples $(x_1,x_3,x_7)$ that satisfy the simplex constraint $x_1 + x_3 + x_7 = 1$.
    Outgoing arcs represent variable assignments to $0$ (dotted) or $1$ (solid).
    Feasible assignments are represented by paths to the truth symbol $\top$.
    For instance, the path marked by green arcs corresponds to $x_1=0, x_3=1, x_7=0$.
}
\label{fig:bdd}
\end{figure}
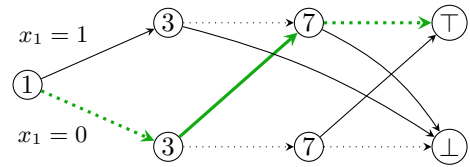

Our definition of BDDs slightly deviates from the common one~\cite{bryant1986graph}, in that we mandate the indices to appear consecutively on each path (even if both outgoing arcs of some node point to the same node).
This property helps us in devising simpler algorithms for min-marginal updates and variable restrictions.
Still, the canonicity property holds, i.e.\ there is a one-to-one correspondence between BDDs and Boolean functions via the following path-property of BDDs.

\begin{prop}[Canonicity]
\label{prop:canonicity}
There is a one-to-one correspondence between Boolean functions $f \colon \{0,1\}^\I \rightarrow \{0,1\}$ and BDDs such that $f(x_{i_1},\ldots,x_{i_k}) = 1$ if, and only if, the corresponding path to $\bddtrue$ in the BDD takes the $x_{i_\ell}$-arc between $v_\ell$ and $v_{\ell+1}$ for $\ell < k$.
\end{prop}

\subsection{Subproblem Representation with BDDs}

In our method we use BDDs to represent indicator functions~$\1_{\X_j}(x)$.
In order to obtain a viable implementation, the BDD size w.r.t.\ the number of variables needs to be moderate.
On the one hand, there exist linear constraints whose corresponding BDD size is exponential irrespective of the variable order~\cite{abio2012new,knuth2011art}. 
On the other hand, there exist polynomial upper bounds on the size of the BDD representation for many commonly encountered constraints.
These include cardinality constraints, constraints with bounded coefficients and many others~\cite{wegener2000branching,knuth2011art}. 
See Figure~\ref{fig:bdd} for an example where $\X$ is defined by a cardinality constraint which has size linear in the number of variables.
For general inequalities there exist lifting techniques utilizing a coefficient splitting that result in multiple BDDs with overall polynomial size~\cite{abio2012new}. 
Hence, the Lagrangean decomposition can be chosen such that the overall size of the BDD representation is polynomial.
In most structured prediction problems and in particular in the ones we solve in our experiments the coefficients in the constraint matrix are bounded by a small constant.
In this case~\citet{knuth2011art} guarantees BDDs of moderate size and we do not need the advanced techniques of~\citet{abio2012new}.

In order to transform linear inequalities into BDDs, we implement the efficient approach of~\citet{Behle2007,behle2008threshold}, which is generalized by \citet{enumerative_branching_serra_2020}.
Our implementation generates a BDD for a given linear inequality in time $\mathcal O(N \log W)$, where $N$ is the number of nodes and $W$ is the maximal width of any layer in the BDD representation.

The representation of $\1_{\X_j}$ as BDDs comes with the advantage that the min-marginals w.r.t.\ $\X_j$ can be computed by dynamic programming (shortest path search) in the associated BDD with appropriate arc costs.
In fact, each forward and backward pass of Algorithm~\ref{alg:min-marginal-averaging} requires traversing all BDDs only once to compute all min-marginals, as intermediate results from previous iterations can be reused.
We detail the update steps of this procedure below.

\subsection{Update Steps for Min-Marginal Computation}

In order to compute min-marginals by shortest path search, we introduce the following BDD arc costs.
\begin{defn}[BDD arc costs]
Let $(V,A)$ be the BDD representing $\1_{\X_j}$.
For every $uv \in A$ we define its arc cost
\begin{equation}
    \theta_{uv} = \begin{cases} 
        \lambda_{\idx(u)}^j,& uv \in A^1 \\
        0,& uv \in A^0\,.
    \end{cases}
\end{equation}
\end{defn}

\begin{algorithm}
    \KwInput{variable index $i \in \I$}
    \For{$v \in V$ with $\idx(v) = i$}{
    $\overrightarrow{m}_{v} = \min\limits_{u: uv \in A} \{\overrightarrow{m}_{u} + \theta_{uv}\}$\;
}
\caption{Forward BDD update step}
\label{alg:forward-step}
\end{algorithm}

\begin{algorithm}
\KwInput{variable index $i \in \I$}
    \For{$v \in V$ with $\idx(v) = i$}{
    $\overleftarrow{m}_{v} = \min \bigg \{ \overleftarrow{m}_{v^+_0} + \theta_{vv^+_0} , \overleftarrow{m}_{v^+_1} + \theta_{vv^+_1} \bigg \}$\;
}
\caption{Backward BDD update step}
\label{alg:backward-step}
\end{algorithm}

In Algorithms~\ref{alg:forward-step} and~\ref{alg:backward-step} we define dynamic programming update steps for computing min-marginals.
For each node $u \in V$ we compute forward messages $\overrightarrow{m}_u$ and backward messages $\overleftarrow{m}_u$, which are the values of the shortest path between $u$ and $r$, respectively $u$ and $\bddtrue$.
The forward messages are computed from first to last variable index using the already computed forward messages of preceding nodes.
Similarly, the backward messages are computed from last to first variable index using the previously computed backward messages of successive nodes.
Backward message values for terminal nodes are fixed to $\overleftarrow{m}_{\bddfalse} = \infty$ and $\overleftarrow{m}_{\bddtrue} = 0$ and the forward message value of the root node $r$ is fixed to $\overrightarrow{m}_{r} = 0$.
Given an arc $uv \in A$, valid forward message $\overrightarrow{m}_u$ and valid backward message $\overleftarrow{m}_{v}$, we define the arc-marginal as
\begin{equation}
    \label{eq:arc-marginal}
    m_{uv} = \overrightarrow{m}_u + \theta_{uv} + \overleftarrow{m}_{v} \\
\end{equation}
and for $i \in \I \cap [n]$, $\beta \in \{0,1\}$ the min-marginals as
\begin{equation}
\label{eq:marginal-aggregation}
    m^{\beta}_{i} = \min_{uv \in A^{\beta} : \idx(u) = i} \{ m_{uv} \},
\end{equation}


\subsection{Incremental BDD Updates}

With the above algorithms for BDD updates, we implement efficient incremental updates for the min-marginal averaging operation in lines~\ref{alg:min-marginal-averaging-mma-update-start}-\ref{alg:min-marginal-averaging-mma-update-end} of Algorithm~\ref{alg:min-marginal-averaging}.
To this end, replace the box on lines~\ref{alg:min-marginal-averaging-mma-update-start}-\ref{alg:min-marginal-averaging-mma-update-end} in the forward pass by
\begin{equation}
\label{eq:min-marginal-averaging-incremental-forward}
    \fbox{\parbox{0.6\linewidth}{Call Algorithm~\ref{alg:forward-step}$(i, \text{BDD}_{\X_j})$; \\[0.5ex]
    Compute min-marginals via~\eqref{eq:marginal-aggregation};}}
\end{equation}
and similarly in the backward pass by
\begin{equation}
\label{eq:min-marginal-averaging-incremental-backward}
    \fbox{\parbox{0.6\linewidth}{Compute min-marginals via~\eqref{eq:marginal-aggregation}; \\[0.5ex] Call Algorithm~\ref{alg:backward-step}$(i, \text{BDD}_{\X_j})$;}}
 \,.
\end{equation}

\begin{prop}
\label{prop:correctness}
After initialization of all messages, the incremental computation in~\eqref{eq:min-marginal-averaging-incremental-forward} and~\eqref{eq:min-marginal-averaging-incremental-backward} produces correct marginals when used in Algorithm~\ref{alg:min-marginal-averaging}.
\end{prop}

With this scheme, we can implement one forward resp.\ backward pass of Algorithm~\ref{alg:min-marginal-averaging} in time proportional to visiting each BDD node exactly once.
In other words, Algorithm~\ref{alg:min-marginal-averaging} has iteration complexity linear in the size of the BDDs that encode the problem decomposition.
For many problems with simple constraints (e.g.\ simplex constraints) the size of BDDs is proportional to the number of variables in the constraint, resulting in complexity linear in the size of the problem description.



\section{Parallelization}
\label{sec:parallelization}
In this section we present a concurrent extension of Algorithm~\ref{alg:min-marginal-averaging} in order to speed up optimization of the dual problem~\eqref{eq:dual-problem}.

The basic idea of our approach is to partition the variable indices into a number of intervals and perform min-marginal averaging for each interval in parallel.
This requires that BDDs which straddle multiple intervals are split into several sub-BDDs, with additional Lagrange dual variables in order to enforce consistency of the representation.
After each forward, resp.\ backward pass the dual variables between sub-BDDs are updated.

More precisely, in this section we assume a partition of $[n] = \{1, \dotsc, n\}$, into $k$ intervals $\mathcal{I}^p = [n_p,n_{p+1})$ with $n_1 = 1$, $n_{k+1} = n+1$ for $p=1,\ldots,k$.
We choose the interval endpoints such that each interval corresponds to a similar number of BDD nodes, see below on how BDD nodes are distributed.

\begin{defn}[Sub-BDD]
\label{def:split-bdd}
    Given a BDD $(V,A)$ we define its \emph{sub-BDDs} $(V_p,A_p)$ for $p=1,\ldots,k$ as follows.
The node set $V_p$ is given by
\begin{align}
V_p & = \{\bddtrue, \bddfalse \} \cup \{v \in V \setsep \idx(v) \in \mathcal{I}^p\} \\
& \quad \cup \{v \in V \setsep \exists vw \in A, \, \idx(w) \in \mathcal{I}^p\}. \nonumber
\end{align}
The set of arcs $A_p$ is given by
\begin{align}
    A_p & = \{ vw \setsep vw \in A, \, \idx(w) \in \mathcal{I}^p \} \\
    & \quad \cup \{ v \bddtrue \setsep vw \in A,\ \idx(v) \in \mathcal{I}^p,\, \idx(w) \notin \mathcal{I}^p\} \nonumber \\
    & \quad \cup \{ v \bddfalse \setsep v \bddfalse \in A \}\,. \nonumber
\end{align}
\end{defn}

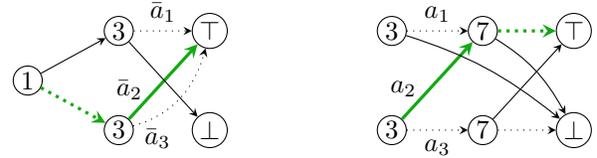
\begin{figure}
\center
\begin{tikzpicture}[xscale=1.1,yscale=1.20,xslant=0.0]
	\tikzstyle{vertex}=[xslant=0.0,draw,circle,inner sep=1pt, minimum width=0pt]         
        
	\node[vertex] (1) at (0,0.5) {$1$};
	\node[vertex] (3a) at (1,0) {$3$};
	\node[vertex] (3b) at (1,1) {$3$};
	\node[vertex] (false) at (2,0) {$\bot$};
	\node[vertex] (true) at (2,1) {$\top$};

	\draw (1) edge[->,>=stealth,dotted,mygreen,very thick] node[above] {} (3a);
	\draw (1) edge[->,>=stealth] node[above] {} (3b);
    \draw (3a) edge[->,>=stealth,dotted,bend right=30] node (3atrue) [below=15,left=0] {$\bar a_3$} (true);
    \draw (3a) edge[->,>=stealth,mygreen,very thick] node (3atrue2) [above=-3,left=4,black] {$\bar a_2$} (true);
    \draw (3b) edge[->,>=stealth,dotted] node (3btrue) [above] {$\bar a_1$} (true);
    \draw (3b) edge[->,>=stealth] node (3bfalse) [above=11, right=-9] {} (false);



	\node[vertex] (3a') at (4,0) {$3$};
	\node[vertex] (3b') at (4,1) {$3$}; 
	\node[vertex] (7a) at (5,0) {$7$};
	\node[vertex] (7b) at (5,1) {$7$};
	\node[vertex] (false') at (6,0) {$\bot$};
	\node[vertex] (true') at (6,1) {$\top$};

    \draw (3a') edge[->,>=stealth,dotted] node (3a7a) [below] {$a_3$} (7a);
    \draw (3a') edge[->,>=stealth,mygreen,very thick] node [above=-3,left=4,black] {$a_2$} (7b);
    \draw (3b') edge[->,>=stealth,dotted] node (3b7b) [above] {$a_1$} (7b);
    \draw (3b') edge[->,>=stealth,bend left=10] node (3bfalse2) [below=-3, left=15] {} (false');
	\draw (7a) edge[->,>=stealth,dotted] (false');
	\draw (7a) edge[->,>=stealth] (true');
	\draw (7b) edge[->,>=stealth,dotted,mygreen,very thick] (true');
	\draw (7b) edge[->,>=stealth,bend left=15] (false');


%
\end{tikzpicture}
\caption{
    Given variable intervals $x_1,x_3 \in \mathcal{I}^p$, $x_7 \in \mathcal{I}^{q}$ with $p < q$, the BDD from Figure~\ref{fig:bdd} is split into two sub-BDDs.
    For consistency between sub-BDDs, active paths must agree on the copy-arc pairs $(\bar a_1, a_1)$, $(\bar a_2, a_2)$ and $(\bar a_3, a_3)$.
    The two paths marked with green arcs correspond to the path marked in Figure~\ref{fig:bdd}.
}
\label{fig:split-bdd}
\end{figure}

The BDD splitting from Definition~\ref{def:split-bdd} is illustrated in Figure~\ref{fig:split-bdd}.
Note that sub-BDDs do not necessarily fulfill condition (i) of Definition~\ref{def:bdd}, as they may have multiple root nodes.
Nonetheless, since they remain acyclic, the update steps detailed in Algorithms~\ref{alg:abstract-forward-step} and \ref{alg:abstract-backward-step} can still be applied in order to compute min-marginals.
Note that arc- and min-marginals now refer to costs w.r.t.\ the sub-BDD, not the original BDD from which those were derived.

\begin{defn}[Copy-arc]
Let $j \in [m]$ be a subproblem index and $(V,A)$ its associated BDD.
Let $uv = a \in A_q$ be an arc of its $q$-th sub-BDD such that $\idx(u) \notin \I^q$ and $v \neq \bot$.
We call the arc $\bar a = u\top \in A_{p}$ for some $p < q$ the \emph{copy-arc} of $a$.
Let~$\C_j = \{ (\bar a, a) \}$ denote the set of all copy-arc pairs of the BDD indexed by $j$.
\end{defn}

For consistency between the sub-BDDs and the original BDD, additional equality constraints between the copy-arcs of the sub-BDDs need to hold (cf.\ Figure~\ref{fig:split-bdd}).
We relax these constraints and reparameterize the associated arc costs with additional dual variables $\mu$.

\begin{defn}
    For any copy-arc pair $(\bar a, a) \in \C_j$ of some BDD $j$ we redefine their arc costs to the values
    \begin{equation}
     \theta^\mu_{\bar a} = \theta_{\bar a} + \mu_{\bar a} \qquad \text{and} \qquad \theta^\mu_a = \theta_a + \mu_a\,.
    \end{equation}
    The additional dual variables $\mu_a, \mu_{\bar a}$ are constrained by
    \begin{equation}
        \mu_a + \mu_{\bar a} \leq 0\,.  
    \end{equation}
\end{defn}

The parallelized min-marginal averaging algorithm is summarized in Algorithm~\ref{alg:min-marginal-averaging-parallel}.
The parallel algorithm works by performing min-marginal averaging in parallel on the intervals $\I^p$. 
After the forward pass, for each copy-arc pair $(\bar a, a) \in \C_j$, the dual variables $\mu_{\bar a}$ are updated with the arc-marginal differences from the forward and backward pass.
This step is possibly dampened by a factor $\gamma \in (0,1]$.
After the backward pass, analogous steps are performed.
Note that the computation of the dual variable updates $\delta^{\rightarrow}$ and $\delta^{\leftarrow}$ in lines~\ref{algl:mma-between-interval-msg-1} and~\ref{algl:mma-between-interval-msg-2} of Algorithm~\ref{alg:min-marginal-averaging-parallel} does not need additional calculations, since the needed arc-marginals are valid after each forward resp.\ backward pass.

\begin{algorithm}
    \caption{Parallel min-marginal averaging}
    \label{alg:min-marginal-averaging-parallel}
Set update step size $\gamma \in (0,1]$\;
Initialize copy-arc updates $\delta^{\rightarrow}_{(\bar a,a)} = \delta^{\leftarrow}_{(\bar a,a)} = 0$\;
    \While{(stopping criterion not met)}{
        \For{$p=1,\ldots,k$ in parallel}{
            Perform forward pass of Algorithm~\ref{alg:min-marginal-averaging}\;
        }

        \For{$p=1,\ldots,k$}{
            \For{all sub-BDDs $j$ in the $p$-th interval}{
                \For{all copy-arc pairs $(\bar a, a) \in \C_j$}{
                    $\delta^{\rightarrow}_{(\bar a, a)} =  m_{\bar a} - \min_{(\bar a, a) \in \C_j} m_{\bar a}$ \label{algl:mma-between-interval-msg-1}\;
        $\mu_{\bar a} \leftarrow \mu_{\bar a} - \gamma \cdot \delta^{\rightarrow}_{(\bar a, a)} + \gamma \cdot \delta^{\leftarrow}_{(\bar a, a)}$ \label{algl:mma-between-interval-msg-update-1}\;
    }
}
        }
        
       \For{$p=1,\ldots,k$ in parallel}{
           Perform backward pass of Algorithm~\ref{alg:min-marginal-averaging}
       }
        
                \For{$p=1,\ldots,k$}{
                    \For{all sub-BDDs $j$ in the $p$-th interval}{
                        \For{all copy-arc pairs $(\bar a, a) \in \C_j$}{
        $\delta^{\leftarrow}_{(\bar a, a)} =  m_{a} - \min_{(\bar a, a) \in \C_j} m_{a}$ \label{algl:mma-between-interval-msg-2}\;
        $\mu_{a} \leftarrow \mu_{a} - \gamma \cdot \delta^{\leftarrow}_{(\bar a, a)} + \gamma \cdot \delta^{\rightarrow}_{(\bar a, a)}$ \label{algl:mma-between-interval-msg-update-2}
    }
    }
    }

    } 
\end{algorithm}

\begin{prop}
\label{prop:correctness-parallel}
    Algorithm~\ref{alg:min-marginal-averaging-parallel} is non-decreasing w.r.t.\ the dual bound.
\end{prop}


\section{Experiments}
\label{sec:experiments}

\begin{table*}[t]
\center
\small
\caption{
		For each dataset the table shows average results obtained with \texttt{Gurobi} and our method \texttt{BDD-MP}.
		 The lower bound denotes the \emph{root LP relaxation} for \texttt{Gurobi} and the value of the Lagrangean relaxation computed by \texttt{BDD-MP}.
		The upper bound denotes the best primal solution value found or \textit{n/a} if for at least one instance no feasible solution was found.
		We report the running times to compute the respective lower and upper bounds within the time limit.
		The best results for each dataset are highlighted in bold.
	}
\begin{tabular}{l|lrrrrrrrrr}
\toprule
\multicolumn{2}{c}{} & \multicolumn{2}{c}{Lower Bound (LB)} & \multicolumn{2}{c}{LB Time [s]} & \multicolumn{2}{c}{Upper Bound (UB)} & \multicolumn{2}{c}{UB Time [s]} \\
\cmidrule(lr){3-4} \cmidrule(lr){5-6} \cmidrule(lr){7-8} \cmidrule(lr){9-10}
\multicolumn{2}{c}{Dataset} & \texttt{Gurobi} & \texttt{BDD-MP} & \texttt{Gurobi} & \texttt{BDD-MP} & \texttt{Gurobi} & \texttt{BDD-MP} & \texttt{Gurobi} & \texttt{BDD-MP} \\
\midrule
\multicolumn{2}{c}{\textit{Cell tracking -- small}} & \textbf{--4.382e06} & --4.387e06 & \textbf{0.42} & 6.85 & \textbf{--4.382e06} & --4.337e06 & \textbf{0.69} & 7.20 \\
\multicolumn{2}{c}{\textit{Cell tracking -- large}} & \textbf{--1.545e08} & --1.549e08 & 262.0 & \textbf{89.2} & \textbf{--1.524e08} & --1.515e08 & 1321.8 & \textbf{127.1} \\
\midrule
\multirow{3}{*}{\textit{GM}}
& \textit{Hotel} & \textbf{--4.293e03} & \textbf{--4.293e03} & \textbf{1.78} & 2.05 & \textbf{--4.293e03} & \textbf{--4.293e03} & \textbf{1.93} & 3.87 \\
& \textit{House} & \textbf{--3.778e03} & \textbf{--3.778e03} & \textbf{3.25} & 4.96 & \textbf{--3.778e03} & \textbf{--3.778e03} & \textbf{3.40} & 6.79 \\
& \textit{Worms} & \textbf{--4.849e04} & --4.878e04 & 273.8 &  \textbf{262.8} & \textbf{--4.842e04} & --4.783e04 & 774.5 & \textbf{264.7} \\
\midrule
\multirow{4}{*}{\textit{MRF}}
& \textit{Color-seg} & \textbf{3.085e08} & \textbf{3.085e08} & 65.0 & \textbf{28.4} & \textbf{3.085e08} & 3.086e08 & 66.5 & \textbf{41.9} \\
& \textit{Color-seg-n4} & \textbf{1.976e04} & 1.964e04 & 443.3 & \textbf{49.2} & 2.846e04 & \textbf{2.179e04} & 537.1 & \textbf{57.5} \\
& \textit{Color-seg-n8} & \textbf{1.973e04} & 1.963e04 & 571.8 & \textbf{97.1} & 2.783e04 & \textbf{2.238e04} & 765.3 & \textbf{120.7} \\
& \textit{Object-seg} & \textbf{3.131e04} & 3.125e04 & 752.6 & \textbf{111.4} & 1.498e05 & \textbf{3.152e04} & 753.1 & \textbf{120.2} \\
\midrule
\multirow{2}{*}{\textit{QAPLIB}}
& \textit{small} & 2.913e06 & \textbf{3.675e06} & 1780.3 & \textbf{176.7} & \textbf{5.186e07} & 5.239e07 & 2167.6 & \textbf{180.7} \\
& \textit{large} & 4.512e04 & \textbf{8.172e06} & 3388.3 & \textbf{2594.5} & \textbf{1.431e08} & 1.452e08 & \textbf{3353.2} & 3357.6 \\
\midrule
\multicolumn{2}{c}{\textit{Discrete tomography}} & \textbf{2.536e02} & 2.394e02 & 321.65 & \textbf{102.81} & \textit{n/a} & \textit{n/a} & \textit{n/a} & \textit{n/a} \\
\bottomrule
\end{tabular}
	\label{tab:results}
\end{table*}


We show competitiveness of our solver with Gurobi~\cite{gurobi}, a leading ILP solver, on a diverse set of structured prediction problems.
More details on the corresponding ILP formulations and experimental results can be found in the appendix.

\myparagraph{Datasets}
We selected a variety of problems from structured prediction and combinatorial optimization to demonstrate the versatility of our solver.
Overall we run experiments on $3115$ instances with 15k--26M variables and 10k--8M constraints. 
The statistics of the datasets in terms of number of variables and constraints are detailed in Table~\ref{tab:statistics} in the appendix.
Our benchmark problems can be categorized as follows.
\begin{description}[leftmargin=0.7cm,noitemsep,topsep=0pt,parsep=0pt,partopsep=0pt]
\setlength{\itemsep}{0pt}
\item[\textnormal{\textit{Cell tracking}}] Small and large cell tracking problems from the study~\cite{haller2020primal}.
\item[\textnormal{\textit{Graph matching (GM)}}] Quadratic assignment problems (often called graph matching in the literature) for correspondence in computer vision~\cite{torresani2008feature} (\texttt{hotel}, \texttt{house}) and developmental biology~\cite{kainmueller2014active} (\texttt{worms}).
\item[\textnormal{\textit{Markov Random Field (MRF)}}] Several datasets from the OpenGM~\cite{Kappes2015} benchmark, containing both small and large instances with varying topologies and number of labels.
\item[\textnormal{\textit{QAPLIB}}] The widely used benchmark dataset for quadratic assignment problems used in the combinatorial optimization community~\cite{QAPLIB}.
\item[\textnormal{\textit{Discrete tomography}}] The synthetic discrete tomography dataset introduced in~\cite{kuske2017novel} consisting of a few thousand instances with a varying number of projections and object densities.
\end{description}

\myparagraph{Algorithms}
Below we specify the algorithms for our empirical comparison. 
For the experiments we set a time limit of 10 minutes for \textit{discrete tomography} and 1 hour for all other instances.
\begin{description}[leftmargin=0.7cm,noitemsep,topsep=0pt,parsep=0pt,partopsep=0pt]
\setlength{\itemsep}{0pt}
	\item[\textnormal{\texttt{Gurobi}}~\textnormal{\cite{gurobi}}] We run the dual simplex method in a single thread to first solve the root LP relaxation and afterwards perform branch-and-bound search.
The dual simplex method was the overall best performing LP algorithm on the considered datasets.
We further disable the presolve routine for better comparability.
In fact, presolve shows little benefit for the considered instances but requires significant time to perform.
	\item[\textnormal{\texttt{BDD-MP}}] We create one BDD per linear inequality. For \textit{cell-tracking-large} and \textit{MRF} we use variable reordering and parallelization with $\gamma = 0.5$.
	 We run Algorithm~\ref{alg:min-marginal-averaging} until a minimal relative objective improvement of $10^{-6}$ and afterwards search for a primal solution with the primal heuristic Algorithm~\ref{alg:primal-heuristic} from the appendix.
		We use the non-smooth optimization except for \textit{discrete tomography}, where we use smoothing parameter $\alpha = 0.01$ (cf.\ Section~\ref{sec:smoothing}).	
\end{description}

\subsection{Results}

In Table~\ref{tab:results} we report averaged upper and lower bounds as well as running times for Gurobi and our method.
Additional convergence plots are provided in the appendix.

\begin{figure*}[!t]
\center
\begin{minipage}{0.33\linewidth}
\begin{tikzpicture}
\footnotesize
\begin{axis}[
	ybar,
	every axis plot/.append style={
          bar width=5pt,
          fill
        },
	width=0.88\columnwidth,
	height=0.65\columnwidth,
	ylabel style={yshift=-1.5em},
	ymin=0,
	ymax=500,
	ytick={400,300,200,100},
	grid style=dotted,
	ymajorgrids=true,
	title={\textit{cell-tracking-large}},
	title style={yshift=-1.7em},
	xtick={1,...,5},
	xtick pos=left,
	xticklabel style = {font=\scriptsize},
	xticklabels={1,2,4,8,16},
	xlabel={\# Threads},
	xlabel style={yshift=0.7ex},
	ylabel={Time [s]},
	legend entries={\texttt{Gurobi}, \texttt{BDD-MP}},
	legend style={
		area legend,
		fill=white,
		draw=none,
		font=\tiny,
		at={(0.52, 0.88)},
		anchor=north east,
		opacity=1,
	},
	legend cell align=left
]
\addplot[fill=red,postaction={pattern=north east lines}] coordinates {(1,262)};
\addplot[fill=blue] coordinates {(1,414) (2,263) (3,174) (4,113) (5,89)};
\end{axis}
\begin{axis}[
	width=0.88\columnwidth,
	height=0.65\columnwidth,
	axis y line*=right,
	ylabel near ticks,
	ylabel style={yshift=0.5em},
	ymin=0,
	ymax=5,
	ytick={1,2,3,4},
	yticklabels={1x,2x,3x,4x},
	xtick={},
	xtick pos=left,
	xticklabels={},
	legend entries={Speedup},
	legend style={
		fill=white,
		draw=none,
		font=\scriptsize,
		at={(0.93, 0.53)},
		anchor=north east,
		text opacity=1,
		fill opacity=0,
	},
]
\addplot[mark=triangle,ultra thick,mygreen] coordinates {(1,1) (2,1.57) (3,2.38) (4,3.66) (5,4.65)};
\end{axis}
\end{tikzpicture}
\end{minipage}
\hfill
\begin{minipage}{0.33\linewidth}
\vspace{0.3em}
\begin{tikzpicture}
\small
\begin{axis}[
	ybar,
	every axis plot/.append style={
          bar width=5pt,
          fill
        },
	width=0.88\columnwidth,
	height=0.65\columnwidth,
	ylabel style={yshift=-1.5em},
	ymin=0,
	ymax=500,
	ytick={100,200,300,400},
	grid style=dotted,
	ymajorgrids=true,
	title={\textit{worms}},
	title style={yshift=-1.5em},
	xtick={1,...,5},
	xtick pos=left,
	xticklabel style = {font=\scriptsize},
	xticklabels={1,2,4,8,16},
	xlabel={\# Threads},
	xlabel style={yshift=0.7ex},
	ylabel={Time [s]},
	legend entries={\texttt{Gurobi}, \texttt{BDD-MP}},
	legend style={
		area legend,
		fill=white,
		draw=none,
		font=\tiny,
		at={(0.45, 0.88)},
		anchor=north east,
		opacity=1,
	},
	legend cell align=left
]
\addplot[fill=red,postaction={pattern=north east lines}] coordinates {(1,274)};
\addplot[fill=blue] coordinates {(1,294) (2,332) (3,290) (4,292) (5,316)};
\end{axis}
\begin{axis}[
	width=0.88\columnwidth,
	height=0.65\columnwidth,
	axis y line*=right,
	ylabel near ticks,
	ylabel style={yshift=0.5em},
	ymin=0,
	ymax=5,
	ytick={1,2,3,4},
	yticklabels={1x,2x,3x,4x},
	xtick={},
	xtick pos=left,
	xticklabels={},
	legend entries={Speedup},
	legend style={
		fill=white,
		draw=none,
		font=\scriptsize,
		at={(0.9, 0.9)},
		anchor=north east,
		text opacity=1,
		fill opacity=0,
	},
]
\addplot[mark=triangle,ultra thick,mygreen] coordinates {(1,1) (2,0.89) (3,1.01) (4,1.01) (5,0.93)};
\end{axis}
\end{tikzpicture}
\end{minipage}
\hfill
\begin{minipage}{0.33\linewidth}
\begin{tikzpicture}
\small
\begin{axis}[
	ybar,
	every axis plot/.append style={
          bar width=5pt,
          fill
        },
	width=0.88\columnwidth,
	height=0.65\columnwidth,
	ylabel style={yshift=-1.5em},
	ymin=0,
	ymax=800,
	ytick={100,300,500,700},
	grid style=dotted,
	ymajorgrids=true,
	title={\textit{object-seg}},
	title style={yshift=-1.7em},
	xtick={1,...,5},
	xtick pos=left,
	xticklabel style = {font=\scriptsize},
	xticklabels={1,2,4,8,16},
	xlabel={\# Threads},
	xlabel style={yshift=0.7ex},
	ylabel={Time [s]},
	legend entries={\texttt{Gurobi}, \texttt{BDD-MP}},
	legend style={
		area legend,
		fill=white,
		draw=none,
		font=\tiny,
		at={(0.5, 0.88)},
		anchor=north east,
		opacity=1,
	},
	legend cell align=left
]
\addplot[fill=red,postaction={pattern=north east lines}] coordinates {(1,753)};
\addplot[fill=blue] coordinates {(1,692) (2,453) (3,287) (4,175) (5,111)};
\end{axis}
\begin{axis}[
	width=0.88\columnwidth,
	height=0.65\columnwidth,
	axis y line*=right,
	ylabel near ticks,
	ylabel style={yshift=0.5em},
	ymin=0,
	ymax=8,
	ytick={1,3,5,7},
	yticklabels={1x,3x,5x,7x},
	xtick={},
	xtick pos=left,
	xticklabels={},
	legend entries={Speedup},
	legend style={
		fill=white,
		draw=none,
		font=\scriptsize,
		at={(0.9, 0.91)},
		anchor=north east,
		text opacity=1,
		fill opacity=0,
	},
]
\addplot[mark=triangle,ultra thick,mygreen] coordinates {(1,1) (2,1.53) (3,2.41) (4,3.95) (5,6.23)};
\end{axis}
\end{tikzpicture}
\end{minipage}
\caption{Running time until convergence (left axes) for Gurobi with dual simplex and our method with variable reordering and 1 to 16 threads.
The right axes show the associated speedup factors of our method.
	}
\label{fig:parallelization}
\end{figure*}
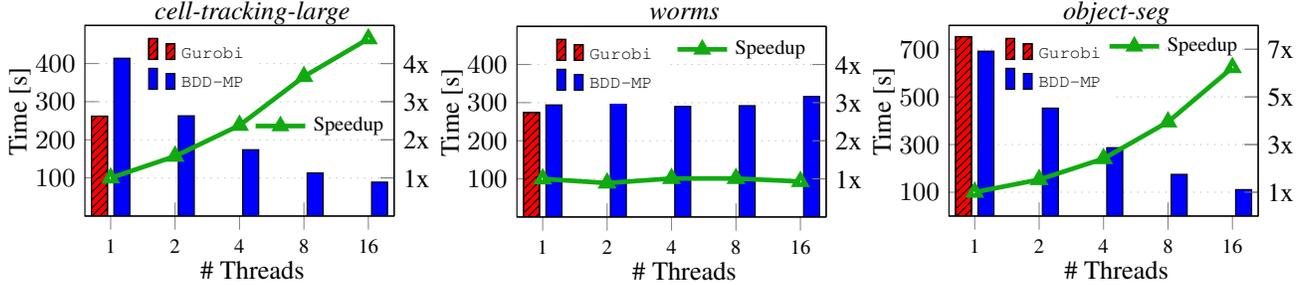

\myparagraph{Lower Bounds}
As can be seen in Table~\ref{tab:results}, for most datasets the lower bounds provided by our method are either equal or slightly weaker compared to those by Gurobi.
The latter is expected, since our solver may converge to suboptimal points of the Lagrangean dual.
Our method spends slightly more time in computing lower bounds for small instances than Gurobi, but is increasingly competitive for larger instances.
In the case of \textit{QAPLIB} both methods provide weak lower bounds while our method performs better than Gurobi.

\myparagraph{Upper Bounds}
The results from Table~\ref{tab:results} indicate that our rather simple primal heuristic achieves small optimality gaps for the considered problems (with the exception of \textit{QAPLIB}).
Generally, our solutions are only slightly worse than those provided by Gurobi.
However, in some cases our method is able to find better primal solutions than Gurobi, which is a consequence of Gurobi's search being restricted by the given time limit.

\myparagraph{Parallelization}
In Figure~\ref{fig:parallelization} we plot the improvement of running time until convergence for our parallelization scheme.
It can be seen that for large instances of \textit{cell tracking} as well as \textit{MRF} the running times decrease considerably with an increasing number of threads.
Naturally, due to significant overhead in the initialization, improvements due to parallelization are only exhibited for large enough problem sizes.
For the denser quadratic assignment problems our current parallelization scheme seems to work less well and requires further research.
In particular the number of iterations until convergence increases, which currently outweighs the fact that each iteration is performed faster.
Generally, large sparse instances of the considered problems benefit from our parallelization scheme.
Note that Gurobi does not make use of multiple threads when solving the relaxation with the dual simplex method.
We provide further evaluation of our parallelization in the appendix.

\myparagraph{Barrier Method}
In contrast to the dual simplex method, the barrier method can be parallelized more effectively.
We tested the performance of Gurobi's barrier method with 16 threads when computing lower bounds for the instances in our study.
For small and sparse problems we found that the dual simplex method is either faster or on par with the barrier method.
For \textit{graph matching} the barrier method ran into numerical problems and could not compute meaningful solutions.
For the larger instances of \textit{QAPLIB} the barrier method could not finish the first iteration within the time limit.
Only for smaller instances of \textit{QAPLIB} the barrier method exhibited superior performance, which we detail in Table~\ref{tab:results-barrier} in the appendix.
Note that starting branch-and-bound search with the dual simplex method given a barrier solution requires an additional crossover step.

\myparagraph{Discussion}
For each of the structured prediction problems in this experimental comparison there are dedicated efficient algorithms that exploit the problem structure, see related work.
They outperform both our method and Gurobi significantly (by at least an order of magnitude), but are not generally applicable and hence unsuitable for new structured prediction problem formulations.

For small instances Gurobi outperforms our solver, while on larger ones we can observe that due to its scalability our solver achieves shorter running times.
The performance of our method in relation to Gurobi is remarkable, as Gurobi and similar leading ILP solvers are highly sophisticated and complex pieces of software that have been developed for decades.
In contrast, our implementation is an academic effort that leaves plenty of room for further improvements in the future.
Our method shows the potential for solving very large scale sparse problems in moderate time due to its scalability and parallelization properties.


\section{Conclusion}
Our contribution serves as a first preliminary step towards more scalable ILP solvers, by providing an approach for solving relaxations of very large scale problems.
Our preliminary results show the potential of our approach in this regard.
However, more work needs to be done in order to make our method competitive, including
(i)~presolving,
(ii)~fine-tuning the initial BDD decomposition,
(iii)~advanced message passing techniques that lead to faster convergence,
(iv)~tightening of the Lagrangean relaxation in terms of additional BDDs,
(v)~integration with branch-and-bound search and
(vi)~parallelization that works uniformly well for sparse and dense problems.

Other open questions include applicability to relaxations of general 0--1 ILPs, i.e.\ the types of problems that can be effectively solved and their optimal BDD representation in terms of size and tightness.

\bibliography{literature.bib}

\begin{thebibliography}{58}
\providecommand{\natexlab}[1]{#1}
\providecommand{\url}[1]{\texttt{#1}}
\expandafter\ifx\csname urlstyle\endcsname\relax
  \providecommand{\doi}[1]{doi: #1}\else
  \providecommand{\doi}{doi: \begingroup \urlstyle{rm}\Url}\fi

\bibitem[Ab{\'\i}o et~al.(2012)Ab{\'\i}o, Nieuwenhuis, Oliveras,
  Rodr{\'\i}guez-Carbonell, and Mayer-Eichberger]{abio2012new}
Ab{\'\i}o, I., Nieuwenhuis, R., Oliveras, A., Rodr{\'\i}guez-Carbonell, E., and
  Mayer-Eichberger, V.
\newblock A new look at bdds for pseudo-boolean constraints.
\newblock \emph{Journal of Artificial Intelligence Research}, 45:\penalty0
  443--480, 2012.

\bibitem[Andersen et~al.(2007)Andersen, Hadzic, Hooker, and
  Tiedemann]{andersen2007constraint}
Andersen, H.~R., Hadzic, T., Hooker, J.~N., and Tiedemann, P.
\newblock A constraint store based on multivalued decision diagrams.
\newblock In \emph{International Conference on Principles and Practice of
  Constraint Programming}, pp.\  118--132. Springer, 2007.

\bibitem[Applegate et~al.(2006)Applegate, Bixby, Chvatal, and
  Cook]{applegate2006traveling}
Applegate, D.~L., Bixby, R.~E., Chvatal, V., and Cook, W.~J.
\newblock \emph{The traveling salesman problem: a computational study}.
\newblock Princeton university press, 2006.

\bibitem[Arora \& Globerson(2013)Arora and Globerson]{arora2013higher}
Arora, C. and Globerson, A.
\newblock Higher order matching for consistent multiple target tracking.
\newblock In \emph{Proceedings of the IEEE International Conference on Computer
  Vision}, pp.\  177--184, 2013.

\bibitem[Bansal et~al.(2004)Bansal, Blum, and Chawla]{bansal2004correlation}
Bansal, N., Blum, A., and Chawla, S.
\newblock Correlation clustering.
\newblock \emph{Machine learning}, 56\penalty0 (1-3):\penalty0 89--113, 2004.

\bibitem[Behle(2007)]{Behle2007}
Behle, M.
\newblock \emph{Binary decision diagrams and integer programming}.
\newblock PhD thesis, Saarland University, 2007.

\bibitem[Behle(2008)]{behle2008threshold}
Behle, M.
\newblock On threshold bdds and the optimal variable ordering problem.
\newblock \emph{Journal of Combinatorial Optimization}, 16\penalty0
  (2):\penalty0 107--118, 2008.

\bibitem[Bergman \& Cire(2016)Bergman and Cire]{BergmanCire2016}
Bergman, D. and Cire, A.~A.
\newblock Decomposition based on decision diagrams.
\newblock In Quimper, C.-G. (ed.), \emph{Integration of AI and OR Techniques in
  Constraint Programming}, pp.\  45--54, Cham, 2016. Springer International
  Publishing.

\bibitem[Bergman \& Cire(2018)Bergman and Cire]{bergman2018discrete}
Bergman, D. and Cire, A.~A.
\newblock Discrete nonlinear optimization by state-space decompositions.
\newblock \emph{Management Science}, 64\penalty0 (10):\penalty0 4700--4720,
  2018.

\bibitem[Bergman et~al.(2015)Bergman, Cire, and van
  Hoeve]{bergman2015lagrangian}
Bergman, D., Cire, A.~A., and van Hoeve, W.-J.
\newblock Lagrangian bounds from decision diagrams.
\newblock \emph{Constraints}, 20\penalty0 (3):\penalty0 346--361, 2015.

\bibitem[Bergman et~al.(2016{\natexlab{a}})Bergman, Cire, Van~Hoeve, and
  Hooker]{bergman2016decision}
Bergman, D., Cire, A.~A., Van~Hoeve, W.-J., and Hooker, J.
\newblock \emph{Decision diagrams for optimization}, volume~1.
\newblock Springer, 2016{\natexlab{a}}.

\bibitem[Bergman et~al.(2016{\natexlab{b}})Bergman, Cire, van Hoeve, and
  Hooker]{bergman2016discrete}
Bergman, D., Cire, A.~A., van Hoeve, W.-J., and Hooker, J.~N.
\newblock Discrete optimization with decision diagrams.
\newblock \emph{INFORMS Journal on Computing}, 28\penalty0 (1):\penalty0
  47--66, 2016{\natexlab{b}}.

\bibitem[Bixby(2012)]{Bixby2012}
Bixby, R.~E.
\newblock A brief history of linear and mixed-integer programming computation.
\newblock \emph{Documenta Mathematica}, Extra Volume: Optimization
  Stories:\penalty0 107--121, 2012.

\bibitem[Bryant(1986)]{bryant1986graph}
Bryant, R.~E.
\newblock Graph-based algorithms for boolean function manipulation.
\newblock \emph{Computers, IEEE Transactions on}, 100\penalty0 (8):\penalty0
  677--691, 1986.

\bibitem[Burkard et~al.(1997)Burkard, Karisch, and Rendl]{QAPLIB}
Burkard, R.~E., Karisch, S.~E., and Rendl, F.
\newblock Qaplib--a quadratic assignment problem library.
\newblock \emph{Journal of Global optimization}, 10\penalty0 (4):\penalty0
  391--403, 1997.

\bibitem[Castro et~al.(2020)Castro, Cire, and Beck]{castro2020mdd}
Castro, M.~P., Cire, A.~A., and Beck, J.~C.
\newblock An mdd-based lagrangian approach to the multicommodity
  pickup-and-delivery tsp.
\newblock \emph{INFORMS Journal on Computing}, 32\penalty0 (2):\penalty0
  263--278, 2020.

\bibitem[{Cplex, IBM ILOG}(2019)]{cplex}
{Cplex, IBM ILOG}.
\newblock Cplex optimization studio 12.10, 2019.

\bibitem[Cuthill \& McKee(1969)Cuthill and McKee]{cuthill1969}
Cuthill, E. and McKee, J.
\newblock Reducing the bandwidth of sparse symmetric matrices.
\newblock In \emph{Proceedings of the 1969 24th National Conference}, ACM
  ’69, pp.\  157–172, New York, NY, USA, 1969. Association for Computing
  Machinery.
\newblock \doi{10.1145/800195.805928}.

\bibitem[Dantzig et~al.(1954)Dantzig, Fulkerson, and
  Johnson]{dantzig1954solution}
Dantzig, G., Fulkerson, R., and Johnson, S.
\newblock Solution of a large-scale traveling-salesman problem.
\newblock \emph{Journal of the operations research society of America},
  2\penalty0 (4):\penalty0 393--410, 1954.

\bibitem[Dlask \& Werner(2020)Dlask and Werner]{dlask2020class}
Dlask, T. and Werner, T.
\newblock A class of linear programs solvable by coordinate-wise minimization.
\newblock \emph{arXiv preprint arXiv:2001.10467}, 2020.

\bibitem[Globerson \& Jaakkola(2008)Globerson and
  Jaakkola]{globerson2008fixing}
Globerson, A. and Jaakkola, T.~S.
\newblock Fixing max-product: Convergent message passing algorithms for {MAP}
  {LP}-relaxations.
\newblock In \emph{Advances in neural information processing systems}, pp.\
  553--560, 2008.

\bibitem[Gondzio \& Sarkissian(2003)Gondzio and
  Sarkissian]{gondzio2003parallel}
Gondzio, J. and Sarkissian, R.
\newblock Parallel interior-point solver for structured linear programs.
\newblock \emph{Mathematical Programming}, 96\penalty0 (3):\penalty0 561--584,
  2003.

\bibitem[Gonz{\'a}lez et~al.(2020)Gonz{\'a}lez, Cire, Lodi, and
  Rousseau]{gonzalez2020bdd}
Gonz{\'a}lez, J.~E., Cire, A.~A., Lodi, A., and Rousseau, L.-M.
\newblock {BDD}-based optimization for the quadratic stable set problem.
\newblock \emph{Discrete Optimization}, pp.\  100610, 2020.

\bibitem[Gonz\'{a}lez et~al.(2020)Gonz\'{a}lez, Cire, Lodi, and
  Rousseau]{gonzalez2020integrated}
Gonz\'{a}lez, J.~E., Cire, A.~A., Lodi, A., and Rousseau, L.-M.
\newblock Integrated integer programming and decision diagram search tree with
  an application to the maximum independent set problem.
\newblock \emph{Constraints}, pp.\  1--24, 2020.

\bibitem[{Gurobi Optimization, LLC}(2020)]{gurobi}
{Gurobi Optimization, LLC}.
\newblock Gurobi optimizer reference manual, 2020.
\newblock URL \url{http://www.gurobi.com}.

\bibitem[Haller et~al.(2020)Haller, Prakash, Hutschenreiter, Pietzsch, Rother,
  Jug, Swoboda, and Savchynskyy]{haller2020primal}
Haller, S., Prakash, M., Hutschenreiter, L., Pietzsch, T., Rother, C., Jug, F.,
  Swoboda, P., and Savchynskyy, B.
\newblock A primal-dual solver for large-scale tracking-by-assignment.
\newblock In \emph{AISTATS}, 2020.

\bibitem[Hooker(2019)]{improved_job_sequence_bounds_hooker_2019}
Hooker, J.~N.
\newblock Improved job sequencing bounds from decision diagrams.
\newblock In Schiex, T. and de~Givry, S. (eds.), \emph{Principles and Practice
  of Constraint Programming}, pp.\  268--283, Cham, 2019. Springer
  International Publishing.
\newblock ISBN 978-3-030-30048-7.

\bibitem[Huangfu \& Hall(2018)Huangfu and Hall]{huangfu2018parallelizing}
Huangfu, Q. and Hall, J. A.~J.
\newblock Parallelizing the dual revised simplex method.
\newblock \emph{Math. Program. Comput.}, 10\penalty0 (1):\penalty0 119--142,
  2018.
\newblock \doi{10.1007/s12532-017-0130-5}.

\bibitem[Jancsary \& Matz(2011)Jancsary and Matz]{jancsary2011convergent}
Jancsary, J. and Matz, G.
\newblock Convergent decomposition solvers for tree-reweighted free energies.
\newblock In \emph{Proceedings of the Fourteenth International Conference on
  Artificial Intelligence and Statistics}, pp.\  388--398, 2011.

\bibitem[Johnson et~al.(2007)Johnson, Malioutov, and
  Willsky]{johnson2007lagrangian}
Johnson, J.~K., Malioutov, D.~M., and Willsky, A.~S.
\newblock Lagrangian relaxation for {MAP} estimation in graphical models.
\newblock \emph{arXiv preprint arXiv:0710.0013}, 2007.

\bibitem[Kainmueller et~al.(2014)Kainmueller, Jug, Rother, and
  Myers]{kainmueller2014active}
Kainmueller, D., Jug, F., Rother, C., and Myers, G.
\newblock Active graph matching for automatic joint segmentation and annotation
  of {C}. elegans.
\newblock In \emph{International Conference on Medical Image Computing and
  Computer-Assisted Intervention}, pp.\  81--88. Springer, 2014.

\bibitem[Kappes et~al.(2015)Kappes, Andres, Hamprecht, Schn{\"{o}}rr, Nowozin,
  Batra, Kim, Kausler, Kr{\"{o}}ger, Lellmann, Komodakis, Savchynskyy, and
  Rother]{Kappes2015}
Kappes, J.~H., Andres, B., Hamprecht, F.~A., Schn{\"{o}}rr, C., Nowozin, S.,
  Batra, D., Kim, S., Kausler, B.~X., Kr{\"{o}}ger, T., Lellmann, J.,
  Komodakis, N., Savchynskyy, B., and Rother, C.
\newblock A comparative study of modern inference techniques for structured
  discrete energy minimization problems.
\newblock \emph{International Journal of Computer Vision}, 115\penalty0
  (2):\penalty0 155--184, 2015.
\newblock \doi{10.1007/s11263-015-0809-x}.

\bibitem[Knuth(2011)]{knuth2011art}
Knuth, D.~E.
\newblock \emph{The art of computer programming, volume 4A: combinatorial
  algorithms, part 1}.
\newblock Pearson Education India, 2011.

\bibitem[Kolmogorov(2006)]{kolmogorov2006convergent}
Kolmogorov, V.
\newblock Convergent tree-reweighted message passing for energy minimization.
\newblock \emph{IEEE transactions on pattern analysis and machine
  intelligence}, 28\penalty0 (10):\penalty0 1568--1583, 2006.

\bibitem[Kolmogorov(2014)]{kolmogorov2014new}
Kolmogorov, V.
\newblock A new look at reweighted message passing.
\newblock \emph{IEEE transactions on pattern analysis and machine
  intelligence}, 37\penalty0 (5):\penalty0 919--930, 2014.

\bibitem[Kuske et~al.(2017)Kuske, Swoboda, and Petra]{kuske2017novel}
Kuske, J., Swoboda, P., and Petra, S.
\newblock A novel convex relaxation for non-binary discrete tomography.
\newblock In \emph{International Conference on Scale Space and Variational
  Methods in Computer Vision}, pp.\  235--246. Springer, 2017.

\bibitem[Lozano et~al.(2018)Lozano, Bergman, and Smith]{lozano2018consistent}
Lozano, L., Bergman, D., and Smith, J.~C.
\newblock On the consistent path problem.
\newblock \emph{Optimization Online e-prints}, 2018.

\bibitem[Luo et~al.(2014)Luo, Xing, Milan, Zhang, Liu, Zhao, and
  Kim]{luo2014multiple}
Luo, W., Xing, J., Milan, A., Zhang, X., Liu, W., Zhao, X., and Kim, T.-K.
\newblock Multiple object tracking: A literature review.
\newblock \emph{arXiv preprint arXiv:1409.7618}, 2014.

\bibitem[Meltzer et~al.(2012)Meltzer, Globerson, and
  Weiss]{meltzer2012convergent}
Meltzer, T., Globerson, A., and Weiss, Y.
\newblock Convergent message passing algorithms-a unifying view.
\newblock \emph{arXiv preprint arXiv:1205.2625}, 2012.

\bibitem[Mittelmann(2017)]{mittelmann2017latest}
Mittelmann, H.~D.
\newblock Latest benchmarks of optimization software.
\newblock In \emph{INFORMS Annual Meeting. Houston, TX}, 2017.

\bibitem[Mittelmann(2020{\natexlab{a}})]{Mittelmann2020}
Mittelmann, H.~D.
\newblock Benchmarking optimization software - a (hi)story.
\newblock \emph{SN Operations Research Forum}, 1, 2020{\natexlab{a}}.
\newblock \doi{10.1007/s43069-020-0002-0}.

\bibitem[Mittelmann(2020{\natexlab{b}})]{MittelmannBenchmark}
Mittelmann, H.~D.
\newblock Benchmarks for optimization software, 2020{\natexlab{b}}.
\newblock URL \url{http://plato.asu.edu/bench.html}.

\bibitem[Papadimitriou(1976)]{papadimitriou1976np}
Papadimitriou, C.~H.
\newblock The np-completeness of the bandwidth minimization problem.
\newblock \emph{Computing}, 16\penalty0 (3):\penalty0 263--270, 1976.

\bibitem[Savchynskyy et~al.(2012)Savchynskyy, Schmidt, Kappes, and
  Schn{\"o}rr]{savchynskyy12efficient}
Savchynskyy, B., Schmidt, S., Kappes, J.~H., and Schn{\"o}rr, C.
\newblock Efficient mrf energy minimization via adaptive diminishing smoothing.
\newblock \emph{UAI. Proceedings}, pp.\  746--755, 2012.
\newblock 1.

\bibitem[Serra(2020)]{enumerative_branching_serra_2020}
Serra, T.
\newblock Enumerative branching with less repetition.
\newblock In \emph{Integration of Constraint Programming, Artificial
  Intelligence, and Operations Research}, pp.\  399--416, 2020.

\bibitem[Swoboda \& Andres(2017)Swoboda and Andres]{swoboda2017message}
Swoboda, P. and Andres, B.
\newblock A message passing algorithm for the minimum cost multicut problem.
\newblock In \emph{Proceedings of the IEEE Conference on Computer Vision and
  Pattern Recognition}, pp.\  1617--1626, 2017.

\bibitem[Swoboda et~al.(2017{\natexlab{a}})Swoboda, Kuske, and
  Savchynskyy]{swoboda2017dual}
Swoboda, P., Kuske, J., and Savchynskyy, B.
\newblock A dual ascent framework for lagrangean decomposition of combinatorial
  problems.
\newblock In \emph{Proceedings of the IEEE Conference on Computer Vision and
  Pattern Recognition}, pp.\  1596--1606, 2017{\natexlab{a}}.

\bibitem[Swoboda et~al.(2017{\natexlab{b}})Swoboda, Rother, Abu~Alhaija,
  Kainmuller, and Savchynskyy]{swoboda2017study}
Swoboda, P., Rother, C., Abu~Alhaija, H., Kainmuller, D., and Savchynskyy, B.
\newblock A study of lagrangean decompositions and dual ascent solvers for
  graph matching.
\newblock In \emph{Proceedings of the IEEE Conference on Computer Vision and
  Pattern Recognition}, pp.\  1607--1616, 2017{\natexlab{b}}.

\bibitem[Swoboda et~al.(2019)Swoboda, Mokarian, Theobalt, Bernard,
  et~al.]{swoboda2019convex}
Swoboda, P., Mokarian, A., Theobalt, C., Bernard, F., et~al.
\newblock A convex relaxation for multi-graph matching.
\newblock In \emph{Proceedings of the IEEE Conference on Computer Vision and
  Pattern Recognition}, pp.\  11156--11165, 2019.

\bibitem[Tjandraatmadja \& van Hoeve(2020)Tjandraatmadja and van
  Hoeve]{tjandraatmadja2020incorporating}
Tjandraatmadja, C. and van Hoeve, W.-J.
\newblock Incorporating bounds from decision diagrams into integer programming.
\newblock \emph{Mathematical Programming Computation}, pp.\  1--32, 2020.

\bibitem[Torresani et~al.(2008)Torresani, Kolmogorov, and
  Rother]{torresani2008feature}
Torresani, L., Kolmogorov, V., and Rother, C.
\newblock Feature correspondence via graph matching: Models and global
  optimization.
\newblock In \emph{European conference on computer vision}, pp.\  596--609.
  Springer, 2008.

\bibitem[Tourani et~al.(2018)Tourani, Shekhovtsov, Rother, and
  Savchynskyy]{tourani2018mplp++}
Tourani, S., Shekhovtsov, A., Rother, C., and Savchynskyy, B.
\newblock Mplp++: Fast, parallel dual block-coordinate ascent for dense
  graphical models.
\newblock In \emph{Proceedings of the European Conference on Computer Vision
  (ECCV)}, pp.\  251--267, 2018.

\bibitem[Tourani et~al.(2020)Tourani, Shekhovtsov, Rother, and
  Savchynskyy]{tourani2020taxonomy}
Tourani, S., Shekhovtsov, A., Rother, C., and Savchynskyy, B.
\newblock Taxonomy of dual block-coordinate ascent methods for discrete energy
  minimization.
\newblock In \emph{AISTATS}, 2020.

\bibitem[Wang \& Koller(2013)Wang and Koller]{wang2013subproblem}
Wang, H. and Koller, D.
\newblock Subproblem-tree calibration: A unified approach to max-product
  message passing.
\newblock In \emph{ICML (2)}, pp.\  190--198, 2013.

\bibitem[Wegener(2000)]{wegener2000branching}
Wegener, I.
\newblock \emph{Branching Programs and Binary Decision Diagrams: Theory and
  Applications}.
\newblock Discrete Mathematics and Applications. Society for Industrial and
  Applied Mathematics, 2000.
\newblock ISBN 9780898719789.
\newblock URL \url{https://books.google.de/books?id=xqqJj42ZoXcC}.

\bibitem[Werner(2007)]{werner2007linear}
Werner, T.
\newblock A linear programming approach to max-sum problem: A review.
\newblock \emph{IEEE transactions on pattern analysis and machine
  intelligence}, 29\penalty0 (7):\penalty0 1165--1179, 2007.

\bibitem[Werner et~al.(2020)Werner, Pr{\r u}{\v s}a, and
  Dlask]{werner2019relative}
Werner, T., Pr{\r u}{\v s}a, D., and Dlask, T.
\newblock Relative interior rule in block-coordinate descent.
\newblock In \emph{Proceedings of the IEEE International Conference on Computer
  Vision}, 2020.
\newblock To appear.

\bibitem[Zhang et~al.(2016)Zhang, Shi, McAuley, Wei, Zhang, and Van
  Den~Hengel]{zhang2016pairwise}
Zhang, Z., Shi, Q., McAuley, J., Wei, W., Zhang, Y., and Van Den~Hengel, A.
\newblock Pairwise matching through max-weight bipartite belief propagation.
\newblock In \emph{Proceedings of the IEEE Conference on Computer Vision and
  Pattern Recognition}, pp.\  1202--1210, 2016.

\end{thebibliography}
\bibliographystyle{icml2021}

\clearpage

\appendix


\section{Lagrangean decomposition}
\label{sec:appendix-decomposition}

\subsection{Derivation of the dual}

We first introduce vectors $y \in \{0,1\}^n$ as well as $x^j \in \X_j$ for $j \in [m]$ and rewrite \eqref{eq:binary-program} redundantly as
\begin{align}
\label{eq:BO-decomposition}
\min_{y,x^1, \dotsc, x^m} \quad & c^\top y \quad
\text{s.t.} \quad  y_{\I_j} = x^j, \; x^j \in \X_j \quad \forall j \in [m].
\end{align}
Now, let $\J_i = \{j \in [m] \mid i \in \I_j\}$ denote the set of variable indices constrained by $\X_j$.
For $i \in [n]$ and $j \in \J_i$ we introduce dual variables $\lambda^j_i$ associated with the equality constraint $y_i = x^j_i$.
For each set of Lagrange variables $\lambda$ we obtain a lower bound to the original problem~\eqref{eq:BO-decomposition} given by
\begin{align}
\min_{y,x^1, \dotsc, x^m} \sum_{i \in [n]} \Big (c_i - \sum_{j \in \J_i} \lambda^j_i \Big ) y_i + \sum_j {\lambda^j}^\top x^j \\
\text{s.t.} \quad  y \in \{0,1\}^n, \quad x^j \in \X_j \quad \forall j \in [m] \nonumber
\end{align}
Optimization above can be decoupled for each $x^j$, $j \in [m]$.
and maximizing over $\lambda$ gives the Lagrangean dual
\begin{align}
\max_\lambda \; \min_{y \in \{0,1\}^n} \sum_{i \in [n]} \Big(c_i - \sum_{j \in \J_i} \lambda^j_i \Big) y_i + \sum_j \min_{x \in \X_j} x^\top \lambda^j.
\end{align}
For simplification we can eliminate $y$ from the formulation by observing that (w.l.o.g.) the maximum is attained for some $\lambda$ that satisfies $\sum_{j \in \J_i} \lambda^j_i = c_i$ for all $i \in [n]$.
Hence, the simplified dual problem reads
\begin{align}
\max_\lambda \quad & \sum_j \min_{x \in \X_j} x^\top \lambda^j \quad
\text{s.t.} \quad \sum_{j \in \J_i} \lambda^j_i = c_i \quad \forall i \in [n]. \tag{D}
\end{align}

\subsection{Min-marginal averaging}

\paragraph{Proof of Proposition \ref{prop:mma-update}}

\begin{proof}
Let $\bar \lambda^j_i = \lambda^j_i - (m^1_{ij} - m^0_{ij})$.
Then 
\begin{align}
E^j(\bar \lambda^j) =
\begin{cases}
E^j(\lambda^j) - (m^1_{ij} - m^0_{ij}) & \text{if } m^1_{ij} - m^0_{ij} < 0 \\
E^j(\lambda^j) & \text{else.}
\end{cases}
\end{align}
Moreover, the min-marginal differences w.r.t.\ $\bar \lambda$ vanish.
Now, let $\bar{\bar \lambda}^j_i = \bar \lambda^j_i + \frac{1}{\abs{\J_i}} \sum_{k \in \J_i} m^1_{ik} - m^0_{ik}$.
Then
\begin{align}
E^j(\bar{\bar \lambda}^j) =
E^j(\bar \lambda^j) + \frac{1}{\abs{\J_i}} \sum_{k \in \J_i} m^1_{ik} - m^0_{ik}
\end{align}
if $\frac{1}{\abs{\J_i}} \sum_{k \in \J_i} m^1_{ik} - m^0_{ik} < 0$ and
\begin{align}
E^j(\bar{\bar \lambda}^j) = E^j(\bar \lambda^j)
\end{align}
otherwise.
Hence, the dual bound increases in total by
\begin{align}
& \sum_{j \in \J_i} E^j(\bar{\bar \lambda}^j) - E^j(\lambda^j) \\
& \quad = - \sum_{\{k \in \J_i \mid m^1_{ik} - m^0_{ik} < 0\}} (m^1_{ik} - m^0_{ik}) \nonumber \\
 & \qquad + \min \left \{ 0, \sum_{k \in \J_i} m^1_{ik} - m^0_{ik} \right \}. \qedhere
\end{align}
\end{proof}


\subsection{Variable order}
\label{sec:variable-order}

The order in which we process the variable indices $i \in [n]$ in Algorithm~\ref{alg:min-marginal-averaging} should facilitate the increase of the dual bound in each iteration.
Therefore, we prefer to process indices $i,i' \in [n]$ consecutively if their updates influence the min-marginals of the associated primal variables $x_i$ and $x_{i'}$, which is the case if there is a subproblem that contains both variables.
A suitable variable order can be obtained by searching for a permutation of the constraint matrix with lowest bandwidth.
The bandwidth of a matrix is the width of the smallest band around the diagonal such that all non-zero entries are contained in it.
Bandwidth-minimization is NP-hard~\cite{papadimitriou1976np}, but fast heuristics such as the Cuthill-McKee algorithm~\cite{cuthill1969} are available.
We run the algorithm on the bipartite variable-constraint adjacency matrix and extract the variable order from the result.


\subsection{Averaging strategy}
\label{sec:averaging-strategy}

The min-marginal averaging update w.r.t.\ $i \in [n]$ defined in \eqref{eq:mma-update} subtracts the min-marginal difference from each corresponding dual variable and distributes the sum of min-marginal differences evenly to all dual variables associated with $i$.
In the case of higher-order graphical models an alternative averaging strategy called Sequential Reweighted Message Passing (SRMP) \cite{kolmogorov2014new} was shown to improve the convergence behavior of the associated DBCA algorithm.
In close analogy to SRMP we suggest the following averaging scheme as an alternative to the default update.
For $i \in [n]$ let
\begin{align}
\J_i^{>} = \{ j \in \J_i \mid \exists i' > i \text{ such that } i' \in \I_j \}
\end{align}
denote the subproblem indices that contain any variable with index greater than $i$, and define $\J_i^{<}$ similarly.
The sets $\J_i^{>}$ and $\J_i^{>}$ are defined here w.r.t.\ the default order on $[n]$ for the sake of simplicity, but can be defined for any other variable order in an analogous way.
The SRMP averaging update distributes the sum of min-marginal differences evenly to all $\lambda^j_i$ for $j \in \J_i^{>}$ in the forward pass.
If $\J_i^{>} = \emptyset$, we fall back to the default averaging scheme by setting $\J_i^{>} = \J_i$.
More precisely, if $\J_i^{>} \neq \emptyset$ during the forward pass, the update \eqref{eq:mma-update} is replaced by
\begin{align}
\lambda^j_i \leftarrow \lambda^j_i - (m^1_{ij} - m^0_{ij}) + \frac{1}{\abs{\J_i^{>}}} \sum_{k \in \J_i} m^1_{ik} - m^0_{ik}
\label{eq:srmp-update}
\end{align}
for $i \in \J_i^{>}$ and
\begin{align}
\lambda^j_i \leftarrow \lambda^j_i - (m^1_{ij} - m^0_{ij})
\end{align}
for $i \in \J_i \setminus \J_i^{>}$.
The backward pass is performed similarly with $\J_i^{<}$ instead of $\J_i^{>}$.


\subsection{Smoothing}
\label{sec:smoothing}
It is well-known that, except in special cases~\cite{dlask2020class}, DBCA can fail to reach the optimum of the relaxation.
Suboptimal stationary points of DBCA algorithms are analyzed in~\cite{werner2019relative}.
One way to attain optima of Lagrangean relaxations with DBCA algorithms is to replace the original non-smooth dual objective with a smooth approximation on which DBCA is guaranteed to find the optimum.
We propose such a smooth approximation below for our Lagrangean decomposition~\eqref{eq:dual-problem} and detail according update rules.
Analogous techniques were used in~\cite{meltzer2012convergent} to smooth the TRWS algorithm~\cite{kolmogorov2006convergent}.

First, we replace the original energy $E^j(\lambda^j)$ through a log-sum-exp based approximation.
For any smoothing parameter $\alpha > 0$ we define
\begin{equation}
    \label{eq:bdd-enery-smoothed}
    E^j_{\alpha}(\lambda^j) = -\alpha \cdot \log\Big(\sum_{x \in \X_j } \exp\Big(\frac{- x^\top \lambda^j}{\alpha}\Big)\Big).
\end{equation}
This results in the smooth Lagrangean dual $\max_{\lambda} \sum_{j \in [m]} E^{j}_{\alpha}(\lambda^j)$. Further, instead of min-marginals~$m_{ij}^\beta$ we define marginal log-sum-exp values as
\begin{equation}
\label{eq:marginal-log-sum-exp}
m_{ij}^{\alpha,\beta} = -\alpha \cdot \log\Big(\sum_{x \in \X_j : x_i = \beta} \exp\Big(\frac{- x^\top \lambda^j}{\alpha}\Big)\Big).
\end{equation}
Finally, the min-marginal averaging operation~\eqref{eq:mma-update} is replaced by
\begin{equation}
\lambda^j_i \leftarrow \lambda^j_i - (m_{ij}^{\alpha,1} - m_{ij}^{\alpha,0}) + \frac{1}{\abs{\J_i}} \sum_{k \in \J_i} m_{ik}^{\alpha,1} - m_{ik}^{\alpha,0}.
\end{equation}
With this change of operations, Algorithm~\ref{alg:min-marginal-averaging} becomes a DBCA method for the smooth approximation.

\begin{prop}[Approximation guarantees]
\label{prop:approx-guarantees}
For any $\alpha > 0$ and $j \in [m]$ it holds that
\begin{equation}
   E^j(\lambda^j) > E^j_{\alpha}(\lambda^j) \geq E^j(\lambda^j) - \alpha \log \: \abs{\X_j}.
\end{equation}
\end{prop}

\begin{proof}
It holds that
\begin{align}
E^j(\lambda^j)
& = \min_{x \in \X_j} x^\top \lambda^j \\
& \geq - \alpha \log \left ( \exp \Big ( - \min_{x \in \X_j} \frac{x^\top \lambda^j}{\alpha} \Big ) \right ) \\
& \geq - \alpha \log \left ( \sum_{x \in \X_j} \exp \Big (-\frac{x^\top \lambda^j}{\alpha} \Big ) \right ) = E_{\alpha}^j(\lambda^j) \\ 
& \geq - \alpha \log \left ( \abs{\X_j} \cdot \exp \Big ( - \min_{x \in \X_j} \frac{x^\top \lambda^j}{\alpha} \Big ) \right ) \\
& = E^j(\lambda^j) - \alpha \log( \abs{\X_j} ). \qedhere
\end{align}
\end{proof}


\section{Primal heuristic}
\subsection{Search strategies}

\begin{algorithm*}
    \KwInput{Open variable indices $\I \subset [n]$, restricted subproblems and indices $(\X_j,\I_j)$, $j \in [m]$, scores $(S_i)_{i \in \I}$}
Find variable $i \in \I$ with maximum score $S_i$\;
(feasible, $\I',(\X'_j,\I'_j)_{j \in [m]}$) = Restriction-Propagation($(\X_j,\I_j)_{j \in [m]}, (i,\beta)$)\;
    \If{feasible and $\I' = \varnothing$}{
\textbf{return} solution\;
    }\uElseIf{feasible}{
feasible = Primal-Heuristic($\I',(\X_j,\I_j)_{j \in [m]},(S_i)_{i \in \I'}$)\;
    }\uElseIf{not feasible}{
(feasible, $\I',(\X'_j,\I'_j)_{j \in [m]}$) = Restriction-Propagation($(\X_j,\I_j)_{j \in [m]}, (i,1-\beta)$)\;
    \If{feasible}{
\textbf{return} Primal-Heuristic($\I',(\X_j,\I_j)_{j \in [m]},(S_i)_{i \in \I'}$)\;
    }\Else{
\textbf{return} \texttt{infeasible}\;
    }
    }
\textbf{output} Partial solution to current subproblem or \texttt{infeasible}\;
\caption{Primal-Heuristic($\I,(\X_j,\I_j)_{j \in [m]},(S_i)_{i \in [n]}$)}
\label{alg:primal-heuristic}
\end{algorithm*}

\begin{algorithm}
    \KwInput{Subproblems $\X_j$, indices $\I_j \subset [n]$, $j \in [m]$, index/value pair to fix $(i,\beta)$.}
    \For{$j \in \J_i$}{
$\X^\beta_j = \{x \in \X_j \setsep x_i = \beta\}$\;
$\F = \{(i',\beta') \in \I_j \times \{0,1\} \mid x \in \X^\beta_j \ \Rightarrow x_{i'} = \beta' \}$\;
$\I_j = \{ i' \in \I_j \setsep \nexists \beta' \text{ s.t. } (i',\beta') \in \F\}$\;
    \For{$(i',\beta') \in \F \backslash \{(x,\beta)\}$}{
Restriction-Propagation($(\X_j,\I_j)_{j \in [m]}, (i',\beta')$)\;
}
}
Set feasible = true $\Leftrightarrow \forall j \in [m]: \X_j \neq \varnothing$\;
    \KwOutput{feasible, restricted subproblems/indices $(\X_j,\I_j)_{j \in [m]}$}
    \caption{\\\hspace*{1em} Restriction-Propagation($(\X_j,\I_j)_{j \in [m]}, (i,\beta)$)}
\label{alg:restriction-propagation}
\end{algorithm}

Another indicator of how suitable a variable/value pair is for fixation is the reduction of the number of feasible solutions when a given variable is fixed to some value.
\begin{defn}[Search space reduction coefficient]
For $i \in [n]$ we define the \emph{search space reduction coefficient} as
\begin{equation}
\label{eq:search-space-reduction-coeff}
 R_i = \sum_{j \in \J_i} \abs{ \{x \in \X_j \mid x_i = 1\} } - \abs{ \{x \in \X_j \mid x_i = 0\} } \,.
\end{equation}
\end{defn} 
The quantity $R_i$ indicates the difference in search space reduction between the fixation $x_i = 1$ and $x_i = 0$ across all associated subproblems.
As an alternative choice for the variable scores $S_i$ that determine the order of variable fixations, we propose $S_i = \sign(R_i) M_i$.
The resulting strategy prefers those variables for which the signs of $R_i$ and $M_i$ agree, thus aligning search space reduction in the individual subproblems with the min-marginal evidence.


\section{Implementation with BDDs}
\label{sec:appendix-bdd}

\paragraph{Proof of Proposition \ref{prop:canonicity}}

\begin{proof}
	See~\cite{bryant1986graph}.
	The required changes due to insertion of BDD nodes with equal outgoing edges do not change the proofs.
\end{proof}

\paragraph{Proof of Proposition \ref{prop:correctness}}

\begin{proof}
First, we note that~\eqref{eq:marginal-aggregation} returns correct marginals for variable $i_\ell$ if we have performed Algorithm~\ref{alg:abstract-backward-step} for variables $i_1,\ldots,i_{\ell-1}$ in that order and Algorithm~\ref{alg:abstract-backward-step} for variables $i_k,\ldots,i_{\ell+1}$ in that order.
	The reason is that Algorithms~\ref{alg:abstract-forward-step} and~\ref{alg:abstract-backward-step} are performing dynamic programming steps for the respective problem, i.e.\ for shortest path with the $(+,\min)$-algebra.
	When processing variable $i$ in the forward pass of Algorithm~\ref{alg:min-marginal-averaging}, forward messages $\overrightarrow{m}$ for variables $i'$, $i' < i$ remain valid, and the same holds true for backward messages $\overleftarrow{m}$ for variables $i' > i$.
	Hence, we only have to update the forward messages for variable $i+1$ to obtain correct marginals for variable $i+1$ via~\eqref{eq:marginal-aggregation}.
	An analogous reasoning holds true for the backward pass.
\end{proof}

\subsection{Abstract BDD update steps}

\begin{algorithm}
    \KwInput{variable index $i \in \I$}
    \For{$v \in V$ with $\idx(v) = i$}{
$\overrightarrow{m}_{v} = \Big (\bigotimes\limits_{u: uv \in A^0} \overrightarrow{m}_{u} \oplus \theta_{uv}  \Big ) \otimes \Big ( \bigotimes\limits_{u: uv \in A^1} (\overrightarrow{m}_{u} \oplus\theta_{uv} ) \Big )$\;
    }
\caption{Abstract forward BDD update step}
\label{alg:abstract-forward-step}
\end{algorithm}

\begin{algorithm}
\KwInput{variable index $i \in \I$}
    \For{$v \in V$ with $\idx(v) = i$}{
\vspace{1.05ex}
$\overleftarrow{m}_{v} = \big(\overleftarrow{m}_{v v^+_0} \oplus \theta_{vv^+_0} \big) \otimes \big(\overleftarrow{m}_{vv^+_1} \oplus \theta_{vv^+_1} \big)$
\vspace{1.1ex}
}
\caption{Abstract backward BDD update step}
\label{alg:abstract-backward-step}
\end{algorithm}

\begin{table*}[!t]
	\setlength{\tabcolsep}{15pt}
    \centering
    \begin{tabular}{cccccccccccccccccccc}
    \toprule
    \multicolumn{5}{c}{Abstract} & \multicolumn{5}{c}{Min-marginal \eqref{eq:min-marginals}} & \multicolumn{5}{c}{Marginal log-sum-exp~\eqref{eq:marginal-log-sum-exp}} & \multicolumn{5}{c}{Solution counts} \\
    \cmidrule(lr){1-5} \cmidrule(lr){6-10} \cmidrule(lr){11-15} \cmidrule(lr){16-20}
\multicolumn{5}{c}{$\big (\oplus, \, \otimes, \, \0, \, \1, \, \theta )$} & \multicolumn{5}{c}{$\big (+, \, \min, \, 0, \, \infty, \, [\lambda^j, 0] \big )$} & \multicolumn{5}{c}{$\big (\cdot, \, +, \, 1, \, 0, \, [\exp(-\lambda^j_i / \alpha), 1] \big )$} & \multicolumn{5}{c}{$(\cdot, \, +, \, 1, \, 0, \, [1, 1])$} \\
    \bottomrule
    \end{tabular}
\caption{
    Symbols in Algorithms~\ref{alg:abstract-forward-step} and~\ref{alg:abstract-backward-step} to compute min-marginal~\eqref{eq:min-marginals}, marginal log-sum-exp~\eqref{eq:marginal-log-sum-exp} and solution counts.}
\label{tab:algebras}
\end{table*}

Similar to the min-marginals \eqref{eq:min-marginals}, we can compute marginal log-sum-exp \eqref{eq:marginal-log-sum-exp} and solution counts (and thus $R_i$ in \eqref{eq:search-space-reduction-coeff}) efficiently with incremental BDD update steps.
To this end, one can employ Algorithms~\ref{alg:abstract-forward-step}--\ref{alg:abstract-backward-step}
that are completely analogous to Algorithms~\ref{alg:forward-step}--\ref{alg:backward-step}  by substitution of the symbols $(\oplus, \, \otimes, \, \0, \, \1, \, \theta^0, \, \theta^1 \big )$
with those listed in Table~\ref{tab:algebras}.

\subsection{Variable fixations}
\label{variable-fixations}

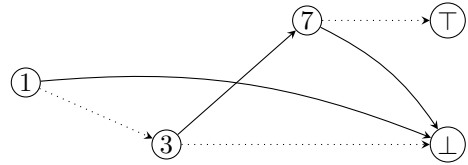
\begin{figure}[!h]
\center
\begin{tikzpicture}[xscale=1.7,yscale=1.5,xslant=0.0]
	\tikzstyle{vertex}=[xslant=0.0,draw,circle,inner sep=1pt, minimum width=0pt]         
        
	\node[vertex] (1) at (0,0.5) {$1$};
	\node[vertex] (3a) at (1,0) {$3$};
	\node[vertex] (7b) at (2,1) {$7$};
	
	\node[vertex] (false) at (3,0) {$\bot$};
	\node[vertex] (true) at (3,1) {$\top$};

	\draw (1) edge[->,>=stealth,dotted] node[above] {} (3a);
	\draw (1) edge[->,>=stealth,bend left=15] node[above] {} (false);
	\draw (3a) edge[->,>=stealth,dotted] (false);
	\draw (3a) edge[->,>=stealth] (7b);
	\draw (7b) edge[->,>=stealth,dotted] (true);
	\draw (7b) edge[->,>=stealth,bend left=15] (false);
		
\end{tikzpicture}
\caption{The BDD obtained from the one depicted in Figure~\ref{fig:bdd} by fixing $x_3 = 1$. Nodes that are no longer reachable from the root or no longer lie on a path to $\top$ have been removed.}
\label{fig:bdd-restricted}
\end{figure}

We implement the variable fixations $x_i = \beta$ in the primal heuristic as manipulations to all BDDs that involve variable~$x_i$.
The manipulations are specified in Algorithms~\ref{alg:variable-fixation} -- \ref{alg:remove-backward}.
In order to fix $x_i = \beta$ in some BDD, Algorithm~\ref{alg:variable-fixation}
redirects all outgoing $(1-\beta)$-arcs from the nodes $v \in V$ with $\idx(v) = i$ to $\bddfalse$.
If that leaves any node unreachable from the root, we recursively remove them with their outgoing arcs by Algorithm~\ref{alg:remove-forward}.
If $\bddtrue$ becomes unreachable, the algorithm detects that the fixation renders the subproblem infeasible.
Similarly, if some node no longer lies on any path towards $\bddtrue$, then we recursively remove it and redirect its incoming arcs by Algorithm~\ref{alg:remove-backward}.
The variable fixation leads to a smaller BDD that represents the restricted feasible set.
See Figure~\ref{fig:bdd-restricted} for an example.

\begin{rem}
The complexity of any sequence of $\leq \abs{\I}$ variable fixations for a BDD is bounded by the number of BDD nodes~$\abs{V}$.
\end{rem}

\begin{algorithm}
    \KwInput{variable index $i \in [n]$, $\beta \in \{0,1\}$}
    \For{$v \in V$ with $\idx(v) = i$}{
Change $\beta$-arc $(v,v^+_\beta)$ to $(v,\bot)$\;
    \If {$v^+_\beta$ has no incoming arcs left}{
    \If{Remove-Forward($v^+_\beta$) = false}{
		\textbf{return} false\;
    }
    }
    \If{both outgoing arcs of $v$ point to $\bot$}{
	Remove-Backward($v$)\;
    }
    }
\textbf{return} true\;
\caption{Variable Fixation}
\label{alg:variable-fixation}
\end{algorithm}

\begin{algorithm}
    \KwInput{BDD node $v \in V$}
    \If{$v = \top$}{
	\textbf{return} false\;
}
Remove $v$ and outgoing arcs $(v,v^+_0)$, $(v,v^+_1)$\;
    \If{$v^+_0$ has no incoming arcs left}{
	 Remove-Forward($v^+_0$)\;
 }
    \If{$v^+_1$ has no incoming arcs left}{
	Remove-Forward($v^+_1$)\;
}
\caption{Remove-Forward}
\label{alg:remove-forward}
\end{algorithm}
\vspace{-4ex}
\begin{algorithm}
    \KwInput{BDD node $v \in V$}
    \For{$(u,v) \in A$}{
	 Replace $(u,v)$ by $(u,\bot)$\;
    \If{both outgoing arcs of $u$ point to $\bot$}{
		Remove-Backward($u$)\;
    }
}
Remove $v$\;
\caption{Remove-Backward}
\label{alg:remove-backward}
\end{algorithm}


\section{Parallelization}
\label{sec:parallelization-appendix}

\paragraph{Proof of Proposition~\ref{prop:correctness-parallel}}
\begin{proof}
    We need to show that the updates of the dual variable $\delta^{\rightarrow}$ and $\delta^{\leftarrow}$ in lines~\ref{algl:mma-between-interval-msg-1} and~\ref{algl:mma-between-interval-msg-2} of Algorithm~\ref{alg:min-marginal-averaging-parallel} results in 
    (i)~feasible dual variables $\mu_{a} + \mu_{\bar a} \leq 0$ for every  copy-arc pair $(\bar a, a)$ and
    (ii)~is non-decreasing in the dual lower bound.
    We will show the statement for the update in line~\ref{algl:mma-between-interval-msg-1}, the update in line~\ref{algl:mma-between-interval-msg-update-2} being analoguous.

    \begin{enumerate}[(i)]
        \item Note that $\delta^{\rightarrow}$ and $\delta^{\leftarrow}$ are non-negative.
            Hence, the first term $-\gamma \cdot \delta^{\rightarrow}_{(\bar a, a)}$ in line~\ref{algl:mma-between-interval-msg-update-1} will decrease $\mu_{\bar a}$ (which does not affect feasibility), while the second term $\gamma \cdot \delta^{\leftarrow}_{(\bar a, a)}$ offsets the changes made in line~\ref{algl:mma-between-interval-msg-update-2}.

        \item 
            The subtraction $-\gamma \cdot \delta^{\rightarrow}_{(\bar a, a)}$ will not decrease the cost of the optimal path in the BDD. First, since $\delta^{\rightarrow}$ is $0$ for the arc the optimal solution takes, the cost of the optimal solution stays the same.
            For the other arcs $\bar a$ that are not in the optimal path, the value $\delta^{\rightarrow}_{(\bar a, a)}$ is the difference of costs of the best path taking the arc $\bar a$ minus taking the optimal path. Hence, the cost update cannot result in previously non-optimal arcs to become optimal ones for $\gamma \leq 1$.
            Hence, the dual lower bound does not decrease after the subtraction.

            Again note that $\delta^{\leftarrow}$ is non-negative. 
            Hence, the second term $\gamma \cdot \delta^{\leftarrow}_{(\bar a, a)}$ in line~\ref{algl:mma-between-interval-msg-update-1} will be non-decreasing in the dual lower bound, since it can only increase the costs.
    \end{enumerate}
\end{proof}


\section{ILP formulations}
\label{sec:ILP-formulations}
Below we detail the ILP formulations of the four problem types considered in the experiments.

\subsection{Markov random fields}

For MRFs, we formulate the associated optimization problem via the local polytope relaxation~\cite{werner2007linear}.
\begin{defn}[MRF formulation]
    Given a graph $G=(V,E)$ with label space $\mathcal{L}_i$ for all $i \in V$, unary potentials $\theta_i \in \R^{\mathcal{L}_i}$ for $i \in V$ and pairwise potentials $\theta_{ij} \in \R^{\mathcal{L}_i \times \mathcal{L}_j}$ for $ij \in E$ we define the feasible set $\Lambda$ as those vectors
    \begin{equation}
            \mu \in \bigotimes_{i \in V} \{0,1\}^{\mathcal{L}_i} \otimes \bigotimes_{ij \in E} \{0,1\}^{\mathcal{L}_i \times \mathcal{L}_j}
    \end{equation}
that satisfy
\begin{equation}
            \begin{array}{cl}
                \sum_{x_i \in \mathcal{L}_i} \mu_i(x_i) = 1 & \forall i \in V, \\
                \sum_{x_i \in \mathcal{L}_i, x_j \in \mathcal{L}_j} \mu_{ij}(x_i,x_j) = 1 & \forall ij \in E, \\
                \sum_{x_j \in \mathcal{L}_j} \mu_{ij}(x_i,x_j) = \mu_{i}(x_i) & \forall ij \in E, x_i \in \mathcal{L}_i, \\
                \sum_{x_i \in \mathcal{L}_i} \mu_{ij}(x_i,x_j) = \mu_{j}(x_j) & \forall ij \in E, x_j \in \mathcal{L}_j.
            \end{array}
    \end{equation}
\end{defn}

The overall 0--1-optimization problem reads
\begin{equation}
    \min_{\mu \in \Lambda} \quad \sum_{i\in V} \la \theta_i, \mu_i \ra + \sum_{ij \in E} \la \theta_{ij}, \mu_{ij} \ra\,.
\end{equation}

\subsection{Graph matching}
The graph matching instances are about matching two sets of points $L$ and $R$. There are both linear costs $c \in \R^{L \times R}$ as well as pairwise costs $d \in \R^{L \times L \times R \times R}$.
We define the feasible set $\Gamma$ as those vectors
\begin{equation}
        \mu \in \{0,1\}^{L \times R}, \quad 
            \nu \in \{0,1\}^{L \times L \times R \times R}
\end{equation}
that satisfy
\begin{equation}
        \begin{array}{cl}
        \mu \mathbbmss{1} \leq \mathbbmss{1}, \\
            \mu^{\top} \mathbbmss{1} \leq \mathbbmss{1}, \\
            \mu_{lr} = \sum_{l'=1}^L \sum_{r'=1}^R \nu_{l l' r r'} & \forall 1 \leq l \leq L, \; 1 \leq r \leq R, \\
            \mu_{lr} = \sum_{l'=1}^L \sum_{r'=1}^R \nu_{l' l r' r} & \forall 1 \leq l \leq L, \; 1 \leq r \leq R.
        \end{array}
\end{equation}
The 0--1-optimization problem is
\begin{equation}
    \min_{(\mu,\nu) \in \Gamma} \; \la c, \mu \ra + \la d, \nu \ra\,.  
\end{equation}
Whenever we have sparse costs $d$ we sparsify our encoding accordingly by leaving out the corresponding variables $\nu$.

\subsection{Cell tracking}
We use the formulation given in~\cite{haller2020primal}.
\begin{defn}[Cell tracking]
Given a set of nodes $V$ corresponding to possible cell detections,
a set of cell transitions $E \in \binom{V}{2}$ and a set of cell divisions $E' \in \binom{V}{3}$,
    we define variables $x_i \in \{0,1\}$, $i \in V$ to correspond to cell detections, $y_{ij} \in \{0,1\}$, $ij \in E$ to cell transitions and $y'_{ijk} \in \{0,1\}$, $ijk \in E'$ to cell divisions.
    Additional conflict sets $C_l \subset V$, $l \in \{1,\ldots,L\}$ for excluding competing cell detection hypotheses are also given.
The feasible set $\C$ is defined as those vectors
\begin{equation}
            x \in \{0,1\}^V, \quad
            y \in \{0,1\}^E, \quad
            y' \in \{0,1\}^{E'}
\end{equation}
that satisfy
\begin{equation}
        \begin{aligned}
            x_i & = \sum\limits_{j: ij \in E} y_{ij} + \sum\limits_{jk: ijk \in E'} y'_{ijk}  & \forall i \in V, \\
            x_j & = \sum\limits_{i: ij \in E} y_{ij} + \sum\limits_{ik: ijk \in E'} y'_{ijk} + \sum\limits_{ik: ikj \in E'} y'_{ikj} & \forall j \in V, \\
            & \sum\limits_{i \in C_l} x_i \leq 1 \quad \forall l \in \{1,\ldots,L\}.
        \end{aligned}
\end{equation}
\end{defn}
Given cell detection costs $\theta_i \in \R$, $i \in V$, cell transition costs $\theta_{ij} \in R$, $ij \in E$ and cell division costs $\theta_{ijk} \in R$, $ijk \in E'$, the 0--1-optimization problem is
\begin{equation}
    \min_{(x,y,y') \in \mathcal{C}} \; \la \theta, (x,y,y')^\top \ra\,.
\end{equation}

\subsection{Discrete tomography}
The discrete tomography problem is encoded as an MRF with additional tomographic projection constraints.
\begin{defn}[Discrete tomography]
    Given a graph $G=(V,E)$ with label space $\mathcal{L}_i = \{0,1,\ldots,k_i\}$ for all $i \in V$, unary potentials $\theta_i \in \R^{\mathcal{L}_i}$ for $i \in V$, pairwise potentials $\theta_{ij} \in \R^{\mathcal{L}_i \times \mathcal{L}_j}$ for $ij \in E$ and tomographic projection constraints
    $\sum_{i \in P_l} x_i = b_l$ for $l \in \{1,\ldots,L\}$, $P_l \subset V$, $b_l \in \N$, the constraint can be summarized as
    \begin{equation}
        \left\{\mu \in \Lambda : \sum_{i \in P_l} \sum_{x_i \in \mathcal{L}_i} x_i \cdot \mu_{i}(x_i) = b_l \quad \forall l \in \{1,\ldots,L\} \right\}\,.
    \end{equation}
\end{defn}


\section{Experiments}

\begin{table}[t]
\center
\small
\begin{tabular}{l|lrrr}
\toprule
\multicolumn{2}{c}{Dataset} & $N$ & Avg $n$ & Avg $m$ \\
\midrule
\multicolumn{2}{c}{\textit{Cell tracking -- small}} & 10 & 22k & 44k \\
\multicolumn{2}{c}{\textit{Cell tracking -- large}} & 5 & 6.0M & 1.3M  \\
\midrule
\multirow{3}{*}{\textit{GM}}
& \textit{Hotel} & 105 & 379k & 52k \\
& \textit{House} & 105 & 379k & 52k \\
& \textit{Worms} & 30 & 1.5M & 166k \\
\midrule
\multirow{4}{*}{\textit{MRF}}
& \textit{Color-seg} & 3 & 2.1M & 8.2M  \\
& \textit{Color-seg-n4} & 9 & 948k & 3.2M  \\
& \textit{Color-seg-n8} & 9 & 1.1M & 6.4M  \\
& \textit{Object-seg} & 5 & 531k & 1.6M \\
\midrule
\multirow{2}{*}{\textit{QAPLIB}}
& \textit{small} & 105 & 399k & 38k \\
& \textit{large} & 29 & 25.8M & 1.1M \\
\midrule
\multicolumn{2}{c}{\textit{Discrete tomography}} & 2700 & 15k & 11k  \\
\bottomrule
\end{tabular}
	\caption{
		For each dataset the table shows the number of instances $N$, average number of variables $n$, average number of constraints~$m$.
	}
	\label{tab:statistics}
\end{table}

\subsection{Additional plots}

\begin{figure*}[t]
\center
\includegraphics[width=0.3\textwidth]{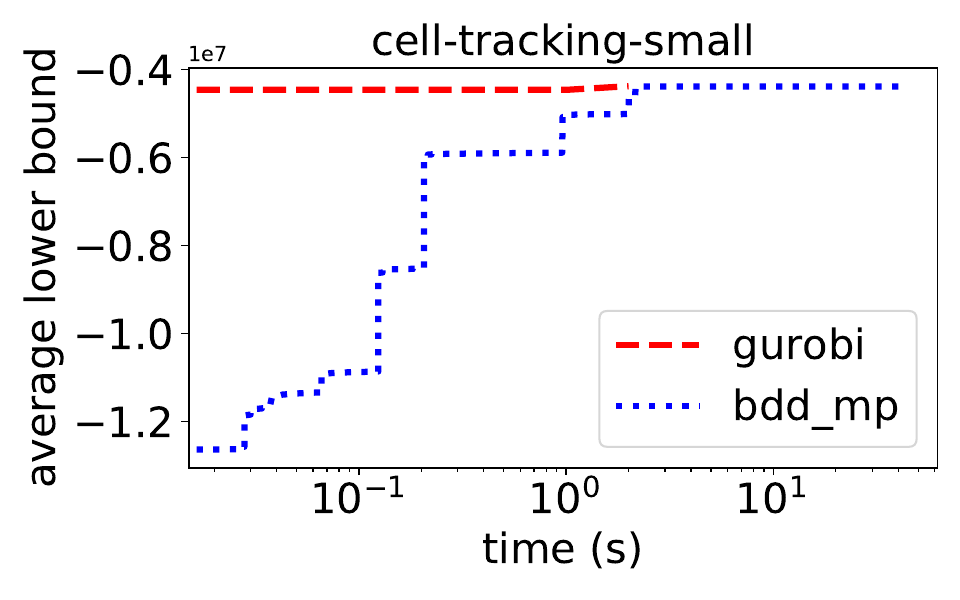}
\includegraphics[width=0.3\textwidth]{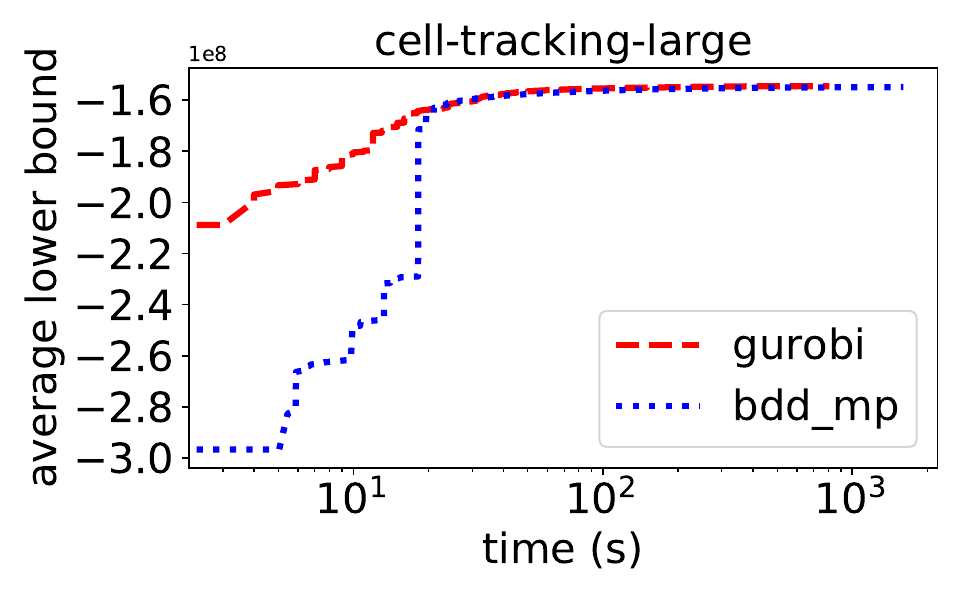}
\includegraphics[width=0.3\textwidth]{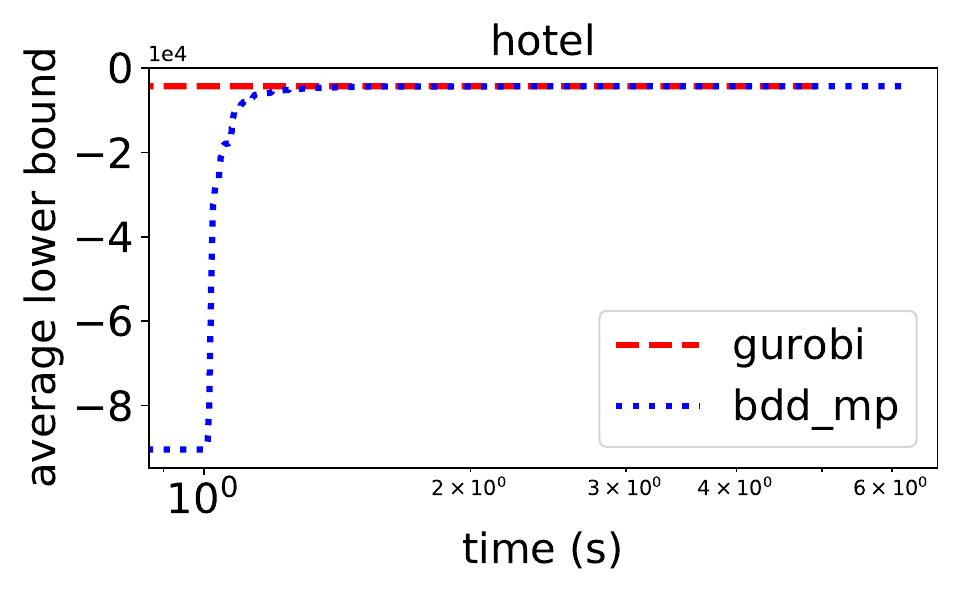}
\hfill
\includegraphics[width=0.3\textwidth]{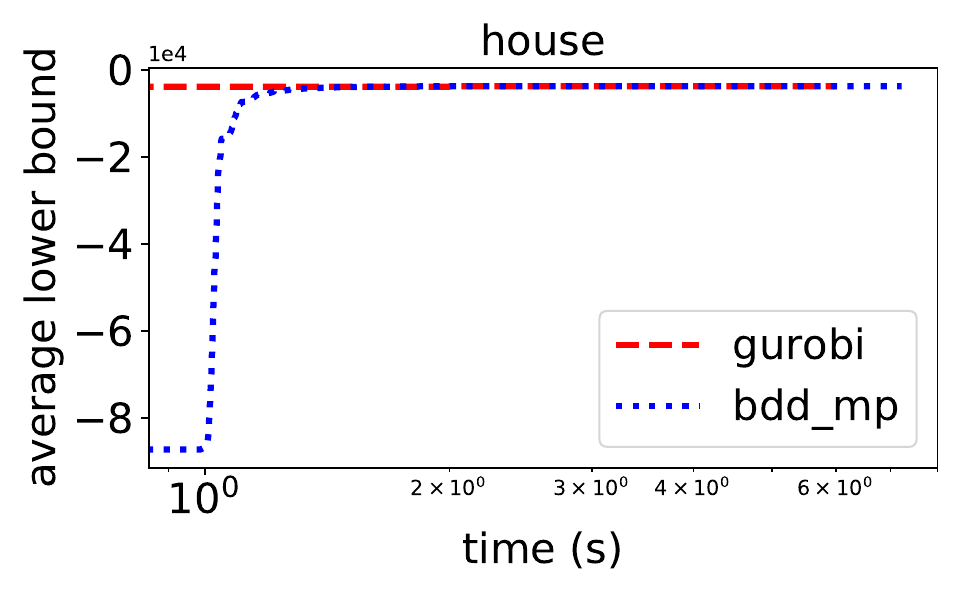}
\includegraphics[width=0.3\textwidth]{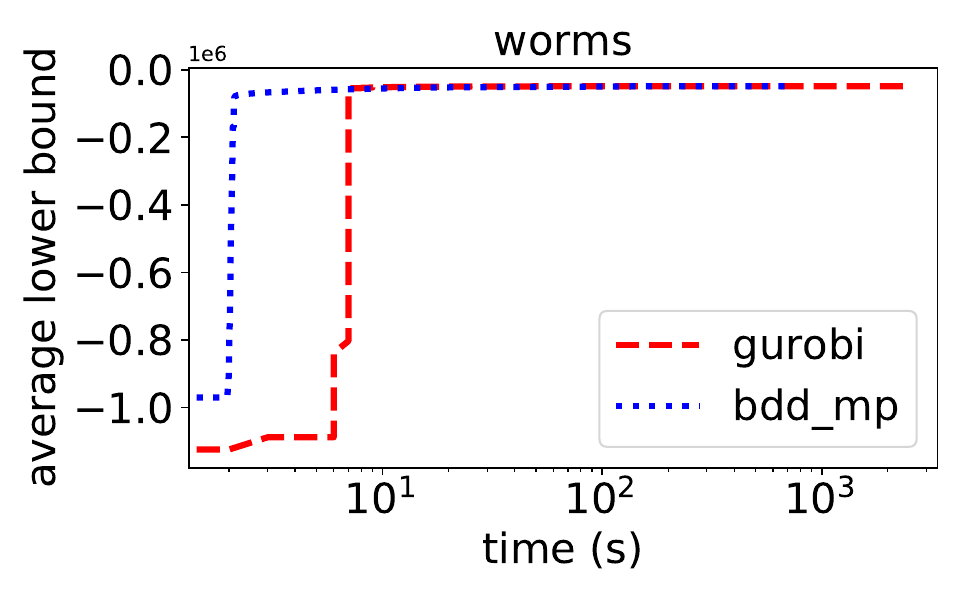}
\includegraphics[width=0.3\textwidth]{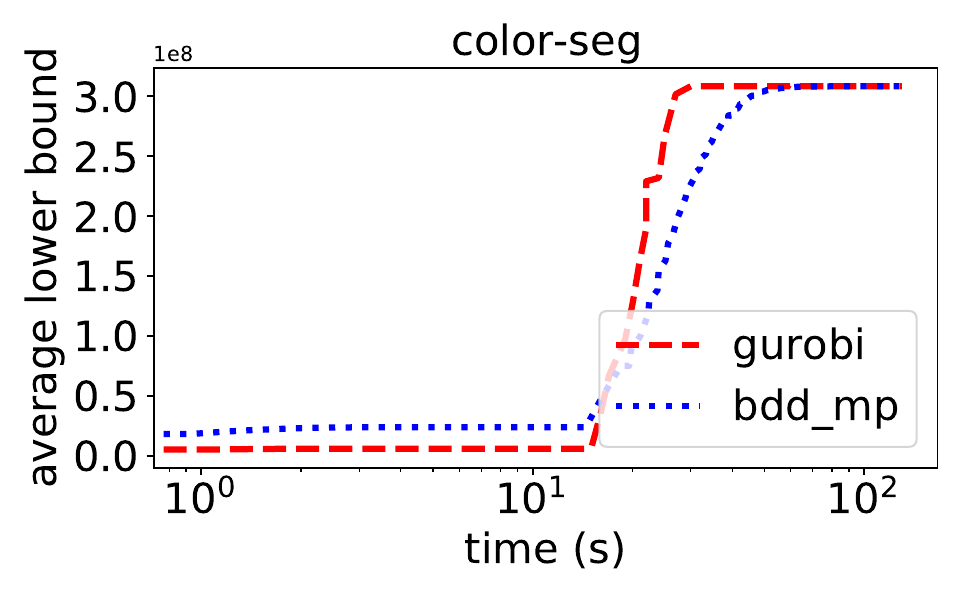}
\hfill
\includegraphics[width=0.3\textwidth]{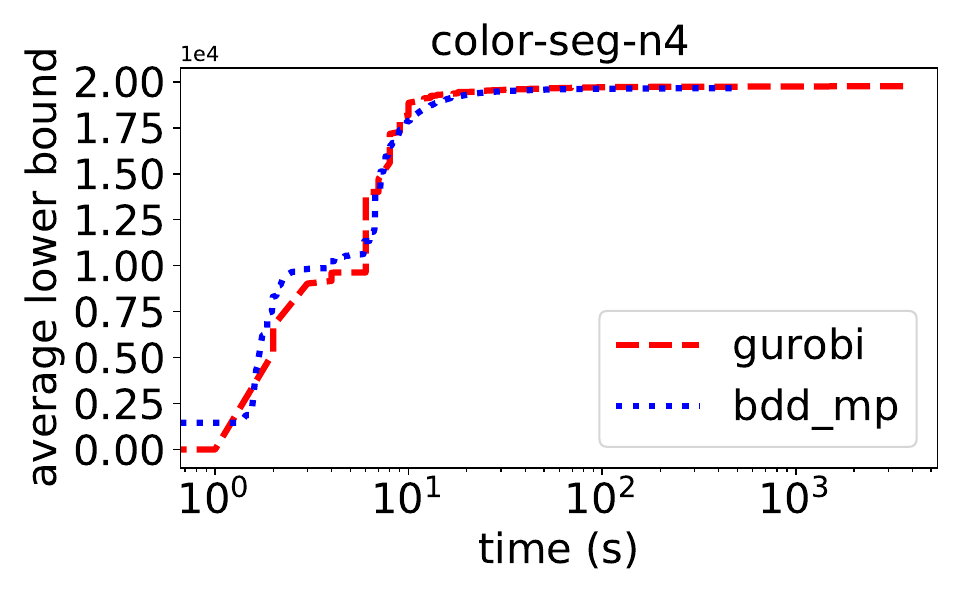}
\includegraphics[width=0.3\textwidth]{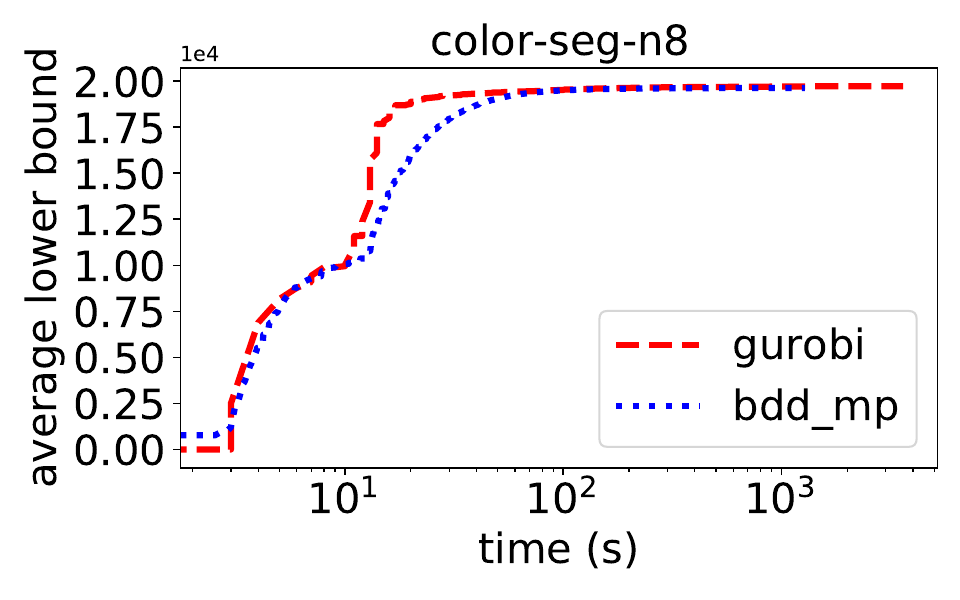}
\includegraphics[width=0.3\textwidth]{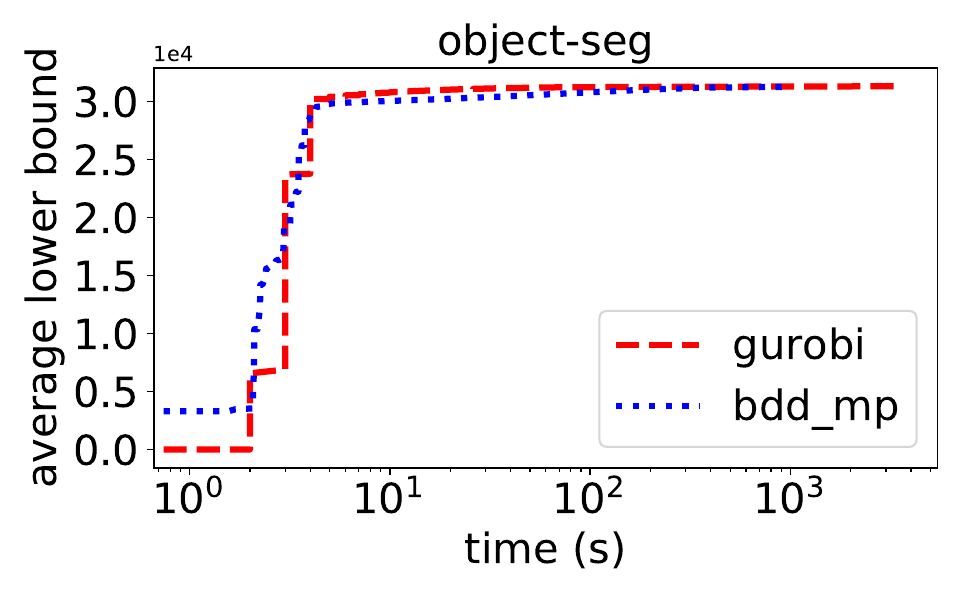}
\hfill
\includegraphics[width=0.3\textwidth]{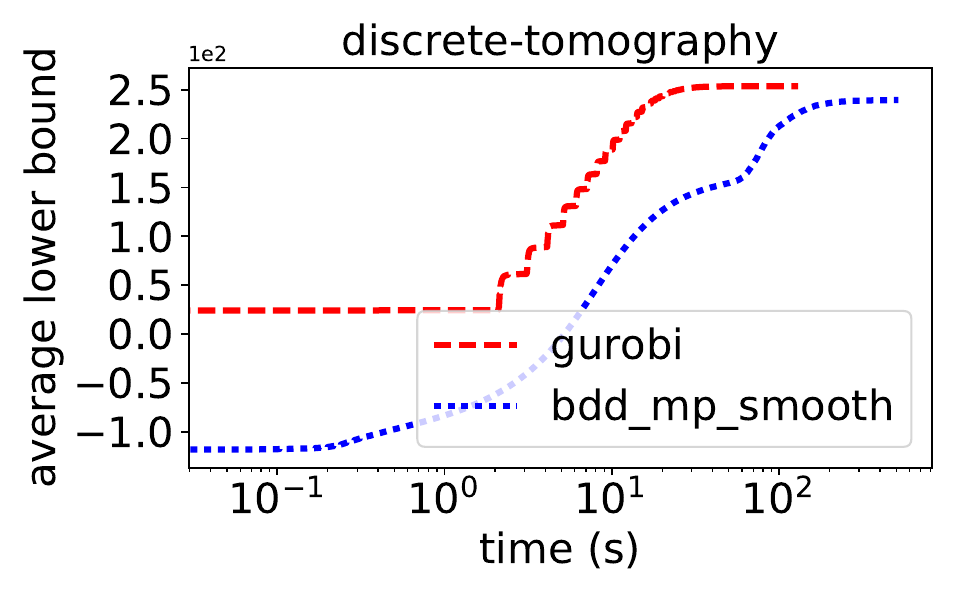}
\includegraphics[width=0.3\textwidth]{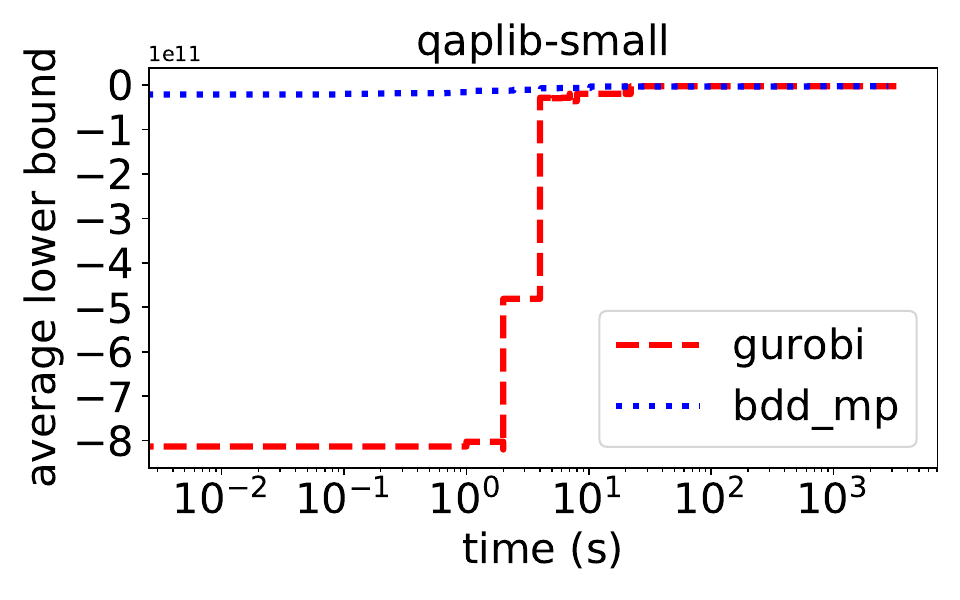}
\includegraphics[width=0.3\textwidth]{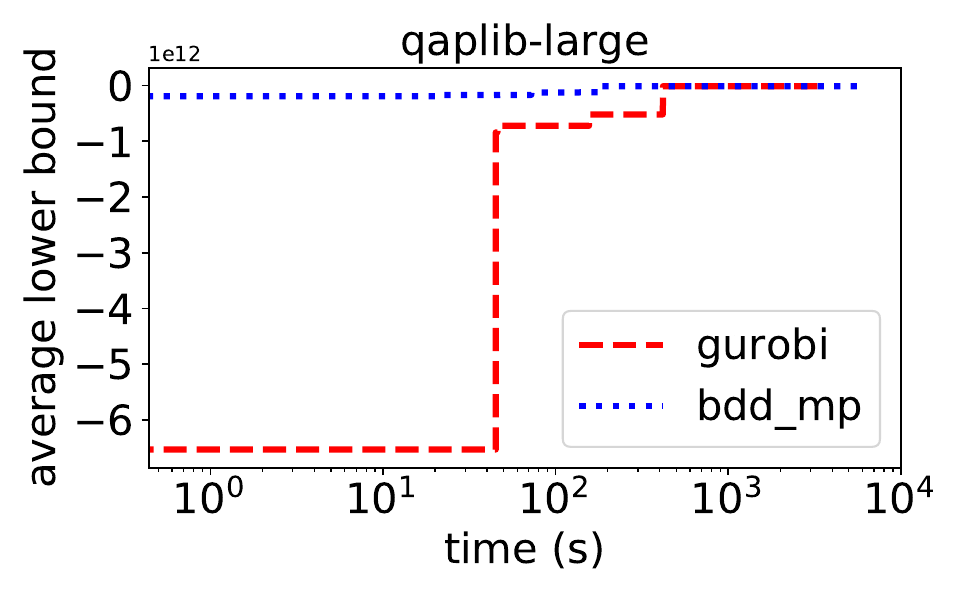}
\caption{Averaged lower bound plots for \texttt{Gurobi} and the basic version of \texttt{BDD-MP} on all datasets.}
\label{fig:lower-bounds-all}
\end{figure*}

\begin{figure*}[!t]
\center
\begin{minipage}{0.33\linewidth}
\begin{tikzpicture}
\small
\begin{axis}[
	ybar,
	every axis plot/.append style={
          bar width=5pt,
          fill
        },
	width=0.88\columnwidth,
	height=0.65\columnwidth,
	ylabel style={yshift=-1.5em},
	ymin=0,
	ymax=100,
	ytick={20,40,60,80},
	grid style=dotted,
	ymajorgrids=true,
	title={\textit{color-seg}},
	title style={yshift=-1.7em},
	xtick={1,...,5},
	xtick pos=left,
	xticklabel style = {font=\scriptsize},
	xticklabels={1,2,4,8,16},
	xlabel={\# Threads},
	xlabel style={yshift=0.7ex},
	ylabel={Time [s]},
	legend entries={\texttt{Gurobi}, \texttt{BDD-MP}},
	legend style={
		area legend,
		fill=white,
		draw=none,
		font=\tiny,
		at={(0.51, 0.88)},
		anchor=north east,
		opacity=1,
	},
	legend cell align=left
]
\addplot[fill=red,postaction={pattern=north east lines}] coordinates {(1,65)};
\addplot[fill=blue] coordinates {(1,82) (2,65) (3,47) (4,41) (5,28)};
\end{axis}
\begin{axis}[
	width=0.88\columnwidth,
	height=0.65\columnwidth,
	axis y line*=right,
	ylabel near ticks,
	ylabel style={yshift=0.5em},
	ymin=0,
	ymax=5,
	ytick={1,2,3,4},
	yticklabels={1x,2x,3x,4x},
	xtick={},
	xtick pos=left,
	xticklabels={},
	legend entries={Speedup},
	legend style={
		fill=white,
		draw=none,
		font=\scriptsize,
		at={(0.92, 0.76)},
		anchor=north east,
		text opacity=1,
		fill opacity=0,
	},
]
\addplot[mark=triangle,ultra thick,mygreen] coordinates {(1,1) (2,1.26) (3,1.74) (4,2.00) (5,2.93)};
\end{axis}
\end{tikzpicture}
\end{minipage}
\hfill
\begin{minipage}{0.33\linewidth}
\begin{tikzpicture}
\small
\begin{axis}[
	ybar,
	every axis plot/.append style={
          bar width=5pt,
          fill
        },
	width=0.88\columnwidth,
	height=0.65\columnwidth,
	ylabel style={yshift=-1.5em},
	ymin=0,
	ymax=700,
	ytick={100,300,500},
	grid style=dotted,
	ymajorgrids=true,
	title={\textit{color-seg-n4}},
	title style={yshift=-1.7em},
	xtick={1,...,5},
	xtick pos=left,
	xticklabel style = {font=\scriptsize},
	xticklabels={1,2,4,8,16},
	xlabel={\# Threads},
	xlabel style={yshift=0.7ex},
	ylabel={Time [s]},
	legend entries={\texttt{Gurobi}, \texttt{BDD-MP}},
	legend style={
		area legend,
		fill=white,
		draw=none,
		font=\tiny,
		at={(0.5, 0.88)},
		anchor=north east,
		opacity=1,
	},
	legend cell align=left
]
\addplot[fill=red,postaction={pattern=north east lines}] coordinates {(1,443)};
\addplot[fill=blue] coordinates {(1,314) (2,194) (3,114) (4,68) (5,49)};
\end{axis}
\begin{axis}[
	width=0.88\columnwidth,
	height=0.65\columnwidth,
	axis y line*=right,
	ylabel near ticks,
	ylabel style={yshift=0.5em},
	ymin=0,
	ymax=7,
	ytick={1,3,5},
	yticklabels={1x,3x,5x},
	xtick={},
	xtick pos=left,
	xticklabels={},
	legend entries={Speedup},
	legend style={
		fill=white,
		draw=none,
		font=\scriptsize,
		at={(0.93, 0.37)},
		anchor=north east,
		text opacity=1,
		fill opacity=0,
	},
]
\addplot[mark=triangle,ultra thick,mygreen] coordinates {(1,1) (2,1.62) (3,2.75) (4,4.62) (5,6.41)};
\end{axis}
\end{tikzpicture}
\end{minipage}
\hfill
\begin{minipage}{0.33\linewidth}
\begin{tikzpicture}
\small
\begin{axis}[
	ybar,
	every axis plot/.append style={
          bar width=5pt,
          fill
        },
	width=0.88\columnwidth,
	height=0.65\columnwidth,
	ylabel style={yshift=-1.5em},
	ymin=0,
	ymax=800,
	ytick={100,300,500,700},
	grid style=dotted,
	ymajorgrids=true,
	title={\textit{color-seg-n8}},
	title style={yshift=-1.7em},
	xtick={1,...,5},
	xtick pos=left,
	xticklabel style = {font=\scriptsize},
	xticklabels={1,2,4,8,16},
	xlabel={\# Threads},
	xlabel style={yshift=0.7ex},
	ylabel={Time [s]},
	legend entries={\texttt{Gurobi}, \texttt{BDD-MP}},
	legend style={
		area legend,
		fill=white,
		draw=none,
		font=\tiny,
		at={(0.53, 0.88)},
		anchor=north east,
		opacity=1,
	},
	legend cell align=left
]
\addplot[fill=red,postaction={pattern=north east lines}] coordinates {(1,572)};
\addplot[fill=blue] coordinates {(1,698) (2,414) (3,272) (4,151) (5,97)};
\end{axis}
\begin{axis}[
	width=0.88\columnwidth,
	height=0.65\columnwidth,
	axis y line*=right,
	ylabel near ticks,
	ylabel style={yshift=0.5em},
	ymin=0,
	ymax=8,
	ytick={1,3,5,7},
	yticklabels={1x,3x,5x,7x},
	xtick={},
	xtick pos=left,
	xticklabels={},
	legend entries={Speedup},
	legend style={
		fill=white,
		draw=none,
		font=\scriptsize,
		at={(0.95, 0.36)},
		anchor=north east,
		text opacity=1,
		fill opacity=0,
	},
]
\addplot[mark=triangle,ultra thick,mygreen] coordinates {(1,1) (2,1.69) (3,2.57) (4,4.62) (5,7.12)};
\end{axis}
\end{tikzpicture}
\end{minipage}
\caption{Running time until convergence (left axes) for Gurobi with dual simplex and our method with variable reordering and 1 to 16 threads.
The right axes show the associated speedup factors of our method.
	}
\label{fig:speedup-appendix}
\end{figure*}
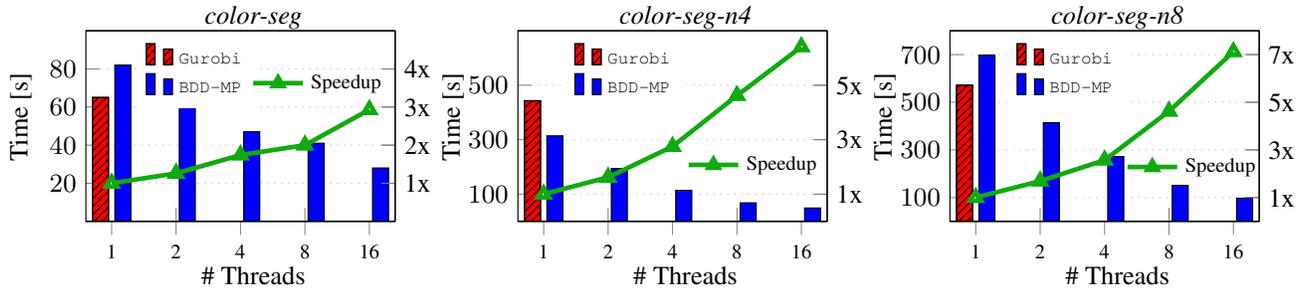

\begin{figure*}[!t]
\center
\includegraphics[width=0.4\textwidth]{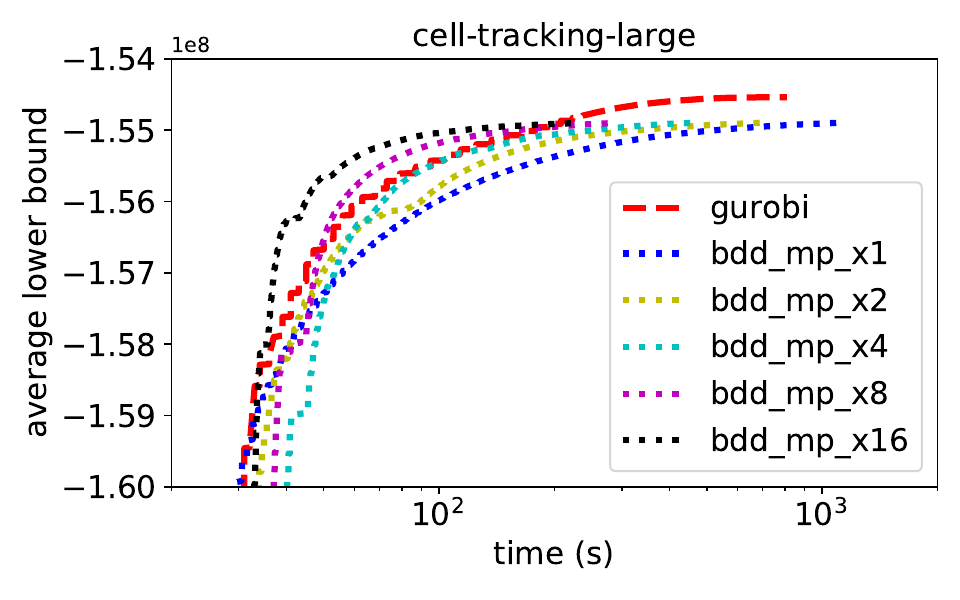}
\hspace{5em}
\includegraphics[width=0.4\textwidth]{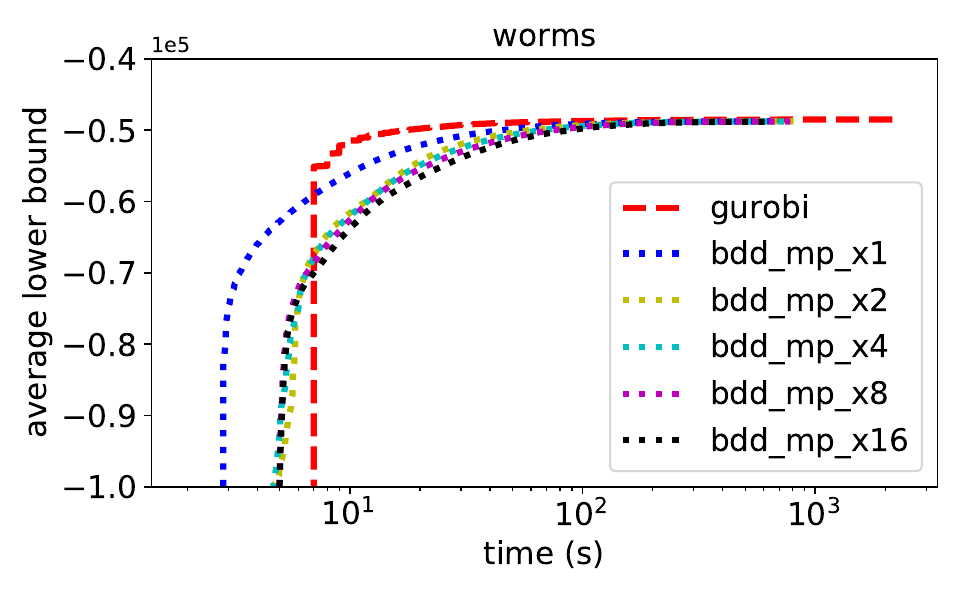}
\includegraphics[width=0.4\textwidth]{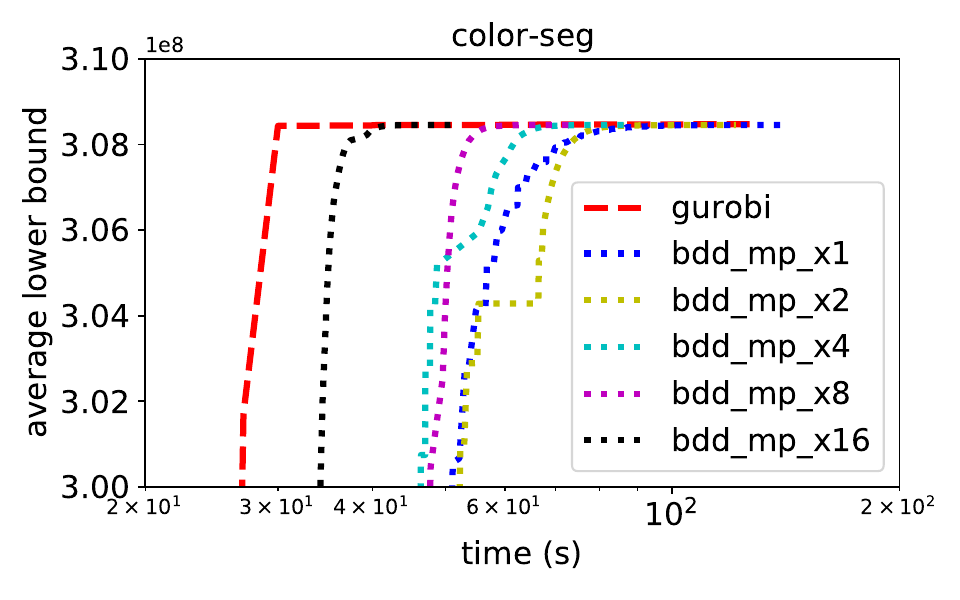}
\hspace{5em}
\includegraphics[width=0.4\textwidth]{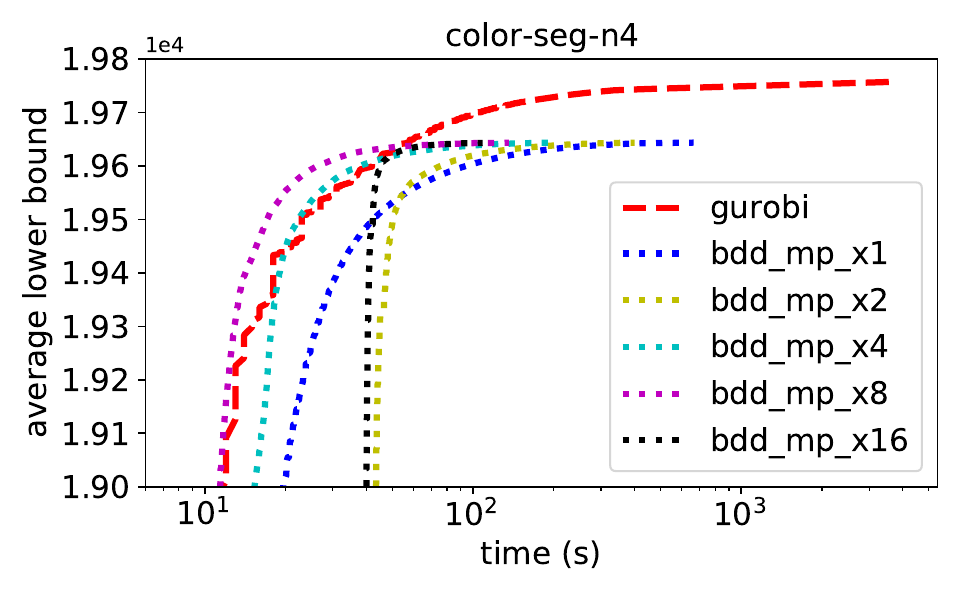}
\includegraphics[width=0.4\textwidth]{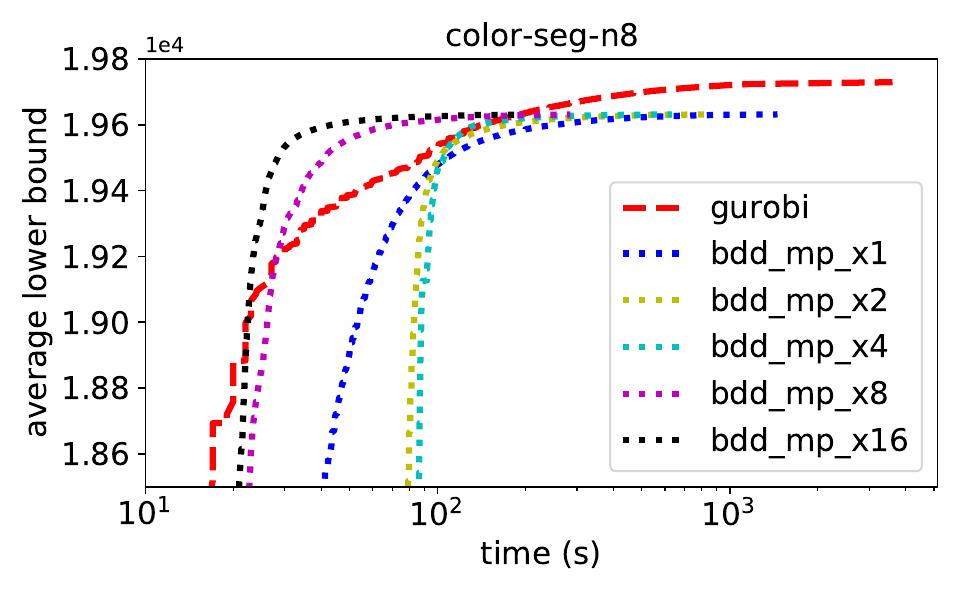}
\hspace{5em}
\includegraphics[width=0.4\textwidth]{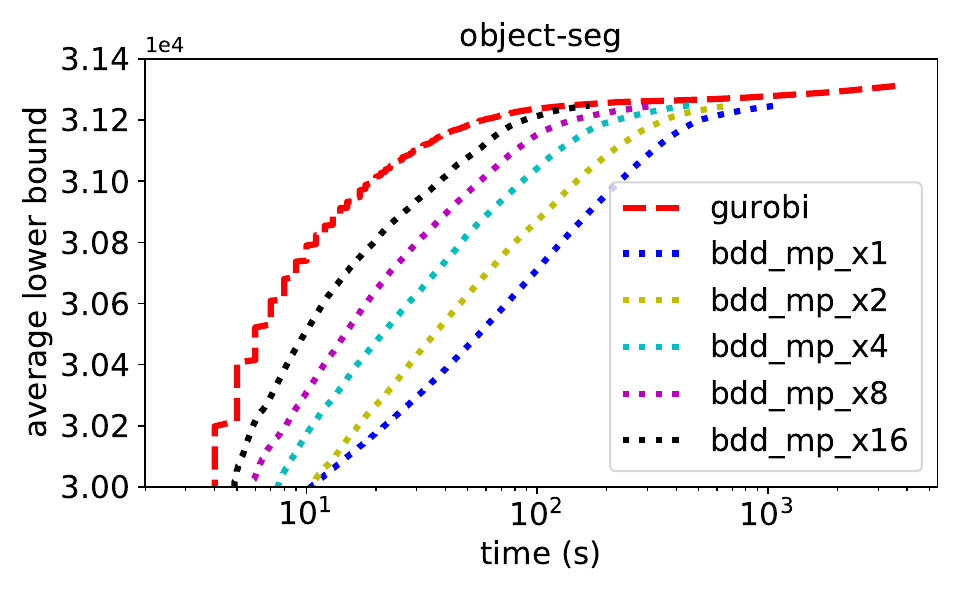}
\caption{
Zoomed-in lower bound convergence plots for Gurobi (dual simplex) and our method with variable reordering and 1 thread (\texttt{x1}) up to 16 threads (\texttt{x16}).
The parallelized versions of our method converge faster after catching up with the initialization overhead.
	}
\label{fig:parallelization-appendix}
\end{figure*}

In Figure \ref{fig:lower-bounds-all}, we show lower bound convergence plots for \texttt{Gurobi} and \texttt{BDD-MP} on all datasets.
In Figure~\ref{fig:speedup-appendix}, we plot the convergence time speedup due to parallelization for the remaining \textit{MRF} instances.
In Figure~\ref{fig:parallelization-appendix}, we show the convergence behavior of lower bounds in relation to the number of threads.

\begin{table*}[t]
\center
\small
\caption{The table reports lower bounds and running times for \textit{QAPLIB small}\footnotemark[1] obtained with Gurobi's dual simplex method (\texttt{Simplex}), its barrier method with 16 threads (\texttt{Barrier}) and our method \texttt{BDD-MP}.
The running times for \texttt{Barrier} do not include any crossover step.
\footnotemark[1]We removed 6 instances that the barrier method could not solve within the time limit in order to enable a comparison.
	}
\begin{tabular}{lrrrrrr}
\toprule
& \multicolumn{3}{c}{Lower Bound (LB)} & \multicolumn{3}{c}{LB Time [s]} \\
\cmidrule(lr){2-4} \cmidrule(lr){5-7}
Dataset & \texttt{Simplex} & \texttt{Barrier} & \texttt{BDD-MP} & \texttt{Simplex} & \texttt{Barrier} & \texttt{BDD-MP} \\
\midrule
\textit{QAPLIB small}\footnotemark[1] & 3.086e06 & \textbf{8.499e06} & 3.554e06 & 1670.2 & 374.6 & \textbf{68.6} \\
\bottomrule
\end{tabular}
	\label{tab:results-barrier}
\end{table*}

\end{document}